\documentclass[12pt]{article}
\usepackage{stmaryrd}
\usepackage{amsmath}
\usepackage{amssymb,amsfonts,amsmath,amsthm,longtable}
\usepackage{epsfig}
\usepackage{pst-poly}     
\usepackage{pst-all}      
\usepackage{amssymb}

\usepackage{cite}

\usepackage{mathrsfs}

\newtheorem{thm}{Theorem}[section]
\newtheorem{cor}{Corollary}[section]

\newtheorem{Def}{Definition}[section]

\makeatletter \@addtoreset{equation}{section} \makeatother
\begin{document}

\title{\bf Fault-Tolerant Path-Embedding of Twisted Hypercube-Like Networks $THLNs$
\thanks {The work was supported  by NNSF of China (No.61472465, 61170303, 61562066) and Natural Science Foundation of Liaoning Province (CN)(No.20170540302).}
\author{Huifeng Zhang, Xirong Xu, Jing Guo, Yuansheng Yang}
\date{\small
{\it  School of Computer Science and Technology\\
\vskip4pt Dalian University of Technology, Dalian, 116024, P.R.China\\}
\vskip4pt
}
}
\maketitle
\begin{abstract}
\noindent The twisted hypercube-like networks($THLNs$) contain several important hypercube variants. This paper is concerned with the fault-tolerant path-embedding of $n$-dimensional($n$-$D$) $THLNs$. Let $G_n$ be an $n$-$D$ $THLN$ and $F$ be a subset of $V(G_n)\cup E(G_n)$ with $|F|\leq n-2$. We show that for arbitrary two different correct vertices $u$ and $v$, there is a faultless path $P_{uv}$ of every length $l$ with $2^{n-1}-1\leq l\leq 2^n-f_v-1-\alpha$, where $\alpha=0$ if vertices $u$ and $v$ form a normal vertex-pair and $\alpha=1$ if vertices $u$ and $v$ form a weak vertex-pair in $G_n-F$($n\geq5$).
\vskip 0.3cm \noindent{\bf Keywords}:\
Multiprocessor interconnection networks, Computer network reliability, Network topology, Hypercubes, Twisted Hypercube-Like Networks $THLNs$, Fault tolerance, Path-embedding.
\end{abstract}
\vskip 0.5cm

\section{Introduction}

The $n$-dimensional hypercube\cite{16},  which possesses many outstanding properties such as recursive structure, relatively small degree, high symmetry, effective routing and broadcasting algorithms\cite{35}, is one of the most efficient, versatile interconnection network and, thus, becomes the preferred topological structure of parallel processing and parallel computing systems\cite{31,33}. Although hypercube networks have many excellent properties, it is well known that they also have inherent shortcomings, such as large diameter. Therefore, many scholars have proposed some hypercube variants, aiming at improving the defects of hypercubes, such as Efe's crossed cubes\cite{6}, Cull's and Larson's Mobius cubes\cite{2}, Hilbers's twisted cubes\cite{15}, Yang's locally twisted cubes\cite{24}. These hypercube variants retain the good properties of hypercubes, but also have many properties superior to hypercubes, such as the diameter of hypercube variants is almost half of the diameter of hypercubes.

Linear arrays (\emph{i.e.} paths), rings (\emph{i.e.} cycles), trees and meshes are the common data structures or foundational interconnection structures used in parallel
computing. The hypercubes and hypercube variants can embed paths\cite{a3, a4, a5}, cycles\cite{a6,a37}, trees\cite{a51,a52}, meshes\cite{a40, a41, a42}. In the process of large-scale Internet operation, it is inevitable that various errors may occur at nodes and edges.
It is significant to find an embedding of a guest graph into a
host graph where all faulty nodes and edges have been removed. This
is called fault-tolerant embedding. Much work has been done on the
fault-tolerant embedding\cite{ 22,36,18,20,3,26,50,128,231,271,101,71,281,131,132,221,601,141,315,316,3128}. A survey paper of Xu and Ma~\cite{36}
lists many results on this topics until 2009.

The hypercube-like networks(short for $HLNs$) are a large class of network topologies \cite{1,18,20}. Among $HLNs$ one
may be identified as a subclass of networks, which in the paper is addressed as the twisted hypercube-like networks (short for $THLNs$), proposed by
Yang\cite{3} in 2011.

 \begin{Def}
\textnormal{\cite{3}
An $n$($n\geq3$)-dimensional (short for $n$-$D$) twisted hypercube-like network (short for $THLN$) is a graph
defined recursively as follows.\\
 \indent(1) A $3$-$D$ $THLN$ is isomorphic to the graph depicted in Fig.1(a).\\
 \indent(2) For $n\geq 4$, an $n$-$D$ $THLN$ $G_n$ is obtained from two vertex-disjoint $(n-1)$-$D$ $THLNs$, denoted by $G^0_{n-1}$ and $G^1_{n-1}$, in this way:
\begin{center}
$V(G_n) = V(G^0_{n-1})\cup V(G^1_{n-1})$,\\
$E(G_n) = E(G^0_{n-1})\cup E(G^1_{n-1})\cup$\\
$ \ \ \ \ \ \ \ \ \ \ \ \ \ \ \ \ \ \ \ \ \ \ \ \ \ \ \{(u,\phi(u)):u\in V(G^0_{n-1})\}$,
\end{center}
where $\phi:V(G^0_{n-1})\rightarrow V(G^1_{n-1})$ is a bijective mapping.
In the following, we will denote this graph $G_n$ as $G_n=\oplus_\phi(G^0_{n-1},G^1_{n-1}$). Fig.1(b) plots a $4D$ $THLN$.}
\end{Def}
\begin{figure}[ht]
\psset{unit=0.8}
\begin{pspicture}(-3,1)(-0.5,2.2)

  \cnode(0.95,2.05){.09}{xl}
  \cnode(0.95,-.15){.09}{vl}
  \cnode(3.15,2.05){.09}{ul}
  \cnode(3.15,-.15){.09}{yl}
  \cnode(0.5,0.95){.09}{wl}
  \cnode(3.6,0.95){.09}{zl}
  \cnode(2.05,2.5){.09}{wla}
  \cnode(2.05,-0.6){.09}{zlb}

  \ncline[linewidth=0.8pt]{xl}{wla}
  \ncline[linewidth=0.8pt]{wl}{xl}
   \ncline[linewidth=0.8pt]{yl}{xl}
   \ncline[linewidth=0.8pt]{zlb}{wla}
   \ncline[linewidth=0.8pt]{ul}{wla}
   \ncline[linewidth=0.8pt]{ul}{zl}
   \ncline[linewidth=0.8pt]{ul}{vl}
   \ncline[linewidth=0.8pt]{wl}{zl}
   \ncline[linewidth=0.8pt]{yl}{zl}
   \ncline[linewidth=0.8pt]{yl}{zlb}
    \ncline[linewidth=0.8pt]{vl}{zlb}
    \ncline[linewidth=0.8pt]{vl}{wl}
    \rput(2.05,-1){\scriptsize(a) $3D$ THLN}
 \end{pspicture}
 \begin{pspicture}(-6.6,1)(-5.9,2.2)

  \cnode(0.95,2.05){.09}{xl}
  \cnode(0.95,-.15){.09}{vl}
  \cnode(3.15,2.05){.09}{ul}
  \cnode(3.15,-.15){.09}{yl}
  \cnode(0.5,0.95){.09}{wl}
  \cnode(3.6,0.95){.09}{zl}
  \cnode(2.05,2.5){.09}{wla}
  \cnode(2.05,-0.6){.09}{zlb}

 \cnode(4.95,2.05){.09}{rxl}
  \cnode(4.95,-.15){.09}{rvl}
  \cnode(7.15,2.05){.09}{rul}
  \cnode(7.15,-.15){.09}{ryl}
  \cnode(4.5,0.95){.09}{rwl}
  \cnode(7.6,0.95){.09}{rzl}
  \cnode(6.05,2.5){.09}{rwla}
  \cnode(6.05,-0.6){.09}{rzlb}

  \ncline[linewidth=0.8pt]{xl}{wla}
  \ncline[linewidth=0.8pt]{wl}{xl}
   \ncline[linewidth=0.8pt]{yl}{xl}
   \ncline[linewidth=0.8pt]{zlb}{wla}
   \ncline[linewidth=0.8pt]{ul}{wla}
   \ncline[linewidth=0.8pt]{ul}{zl}
   \ncline[linewidth=0.8pt]{ul}{vl}
   \ncline[linewidth=0.8pt]{wl}{zl}
   \ncline[linewidth=0.8pt]{yl}{zl}
   \ncline[linewidth=0.8pt]{yl}{zlb}
    \ncline[linewidth=0.8pt]{vl}{zlb}
    \ncline[linewidth=0.8pt]{vl}{wl}

    \ncline[linewidth=0.8pt]{rxl}{rwla}
  \ncline[linewidth=0.8pt]{rwl}{rxl}
   \ncline[linewidth=0.8pt]{ryl}{rxl}
   \ncline[linewidth=0.8pt]{rzlb}{rwla}
   \ncline[linewidth=0.8pt]{rul}{rwla}
   \ncline[linewidth=0.8pt]{rul}{rzl}
   \ncline[linewidth=0.8pt]{rul}{rvl}
   \ncline[linewidth=0.8pt]{rwl}{rzl}
   \ncline[linewidth=0.8pt]{ryl}{rzl}
   \ncline[linewidth=0.8pt]{ryl}{rzlb}
    \ncline[linewidth=0.8pt]{rvl}{rzlb}
    \ncline[linewidth=0.8pt]{rvl}{rwl}

 \ncline[linewidth=.4pt]{wla}{rwla}
  \ncline[linewidth=.4pt]{ul}{rxl}
   \nccurve[angleA=21,angleB=158,linewidth=.4pt]{zl}{rzl}
   \ncline[linewidth=.4pt]{yl}{rul}
   \ncline[linewidth=.4pt]{zlb}{rzlb}
   \nccurve[angleA=21,angleB=158,linewidth=.4pt]{vl}{ryl}
   \nccurve[angleA=21,angleB=158,linewidth=.4pt]{wl}{rwl}
   \ncline[linewidth=.4pt]{xl}{rvl}
 \rput(4.05,-1){\scriptsize (b) $4D$ THLN}
 \end{pspicture}
 \vskip40pt
\caption{
\label{f4}                                    
\footnotesize Examples of $3D$ $THLN$ and $4D$ $THLN$}
\end{figure}
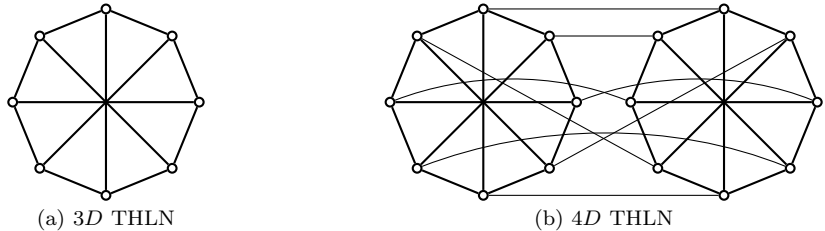

Specifically, the fore-mentioned hypercube variant networks are all $THLNs$. In 2005, Park \emph{et al}.\cite{18} demonstrated that all $n$-$D$ $THLNs$ are Hamiltonian with at most
$n-2$ faulty elements and Hamiltonian connected with at most $n-3$ faulty elements.  Furthermore, Zhang \emph{et al}.\cite{128} improved the upper bound of fault tolerant Hamiltonian connectivity to $n-2$ excepting only a pair of vertices and gave the definitions of weak vertex-pair and normal vertex-pair as follows.

\begin{Def}$^{[15]}$
\textnormal{Let $F\subset V(G_n)\cup E(G_n)$ with $|F|=n-2$. If $G_n-F$ contains a vertex $w$ such that $N_{G_n-F}(w)=\{w_1,w_2\}$, then $w$ is called as a weak 2-degree vertex and $(w_1,w_2)$  is called as a $w$-weak vertex pair(short for weak vertex pair ).}
\end{Def}
If $F=\{a, b\}$, for instance, then $w$ is a
weak 2-degree vertex and $(w_1,w_2)$ is a weak
vertex-pair in $G_4-F$(See Figure 2).
 \vskip10pt
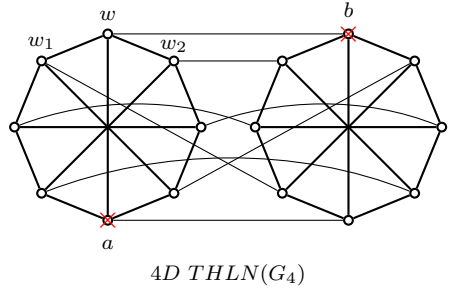
\begin{figure}[ht]
\psset{unit=0.8}
\hskip70pt
 \begin{pspicture}(-2.6,1)(-1.9,2.2)

  \cnode(0.95,2.05){.09}{xl}
  \cnode(0.95,-.15){.09}{vl}
  \cnode(3.15,2.05){.09}{ul}
  \cnode(3.15,-.15){.09}{yl}
  \cnode(0.5,0.95){.09}{wl}
  \cnode(3.6,0.95){.09}{zl}
  \cnode(2.05,2.5){.09}{wla}
  \cnode(2.05,-0.6){.09}{zlb}

  \cnode(4.95,2.05){.09}{rxl}
  \cnode(4.95,-.15){.09}{rvl}
  \cnode(7.15,2.05){.09}{rul}
  \cnode(7.15,-.15){.09}{ryl}
  \cnode(4.5,0.95){.09}{rwl}
  \cnode(7.6,0.95){.09}{rzl}
  \cnode(6.05,2.5){.09}{rwla}
  \cnode(6.05,-0.6){.09}{rzlb}

  \ncline[linewidth=0.8pt]{xl}{wla}
  \ncline[linewidth=0.8pt]{wl}{xl}
   \ncline[linewidth=0.8pt]{yl}{xl}
   \ncline[linewidth=0.8pt]{zlb}{wla}
   \ncline[linewidth=0.8pt]{ul}{wla}
   \ncline[linewidth=0.8pt]{ul}{zl}
   \ncline[linewidth=0.8pt]{ul}{vl}
   \ncline[linewidth=0.8pt]{wl}{zl}
   \ncline[linewidth=0.8pt]{yl}{zl}
   \ncline[linewidth=0.8pt]{yl}{zlb}
    \ncline[linewidth=0.8pt]{vl}{zlb}
    \ncline[linewidth=0.8pt]{vl}{wl}

    \ncline[linewidth=0.8pt]{rxl}{rwla}
   \ncline[linewidth=0.8pt]{rwl}{rxl}
   \ncline[linewidth=0.8pt]{ryl}{rxl}
   \ncline[linewidth=0.8pt]{rzlb}{rwla}
   \ncline[linewidth=0.8pt]{rul}{rwla}
   \ncline[linewidth=0.8pt]{rul}{rzl}
   \ncline[linewidth=0.8pt]{rul}{rvl}
   \ncline[linewidth=0.8pt]{rwl}{rzl}
   \ncline[linewidth=0.8pt]{ryl}{rzl}
   \ncline[linewidth=0.8pt]{ryl}{rzlb}
    \ncline[linewidth=0.8pt]{rvl}{rzlb}
    \ncline[linewidth=0.8pt]{rvl}{rwl}

   \ncline[linewidth=.4pt]{wla}{rwla}
   \ncline[linewidth=.4pt]{ul}{rxl}
   \nccurve[angleA=21,angleB=158,linewidth=.4pt]{zl}{rzl}
   \ncline[linewidth=.4pt]{yl}{rul}
   \ncline[linewidth=.4pt]{zlb}{rzlb}
   \nccurve[angleA=21,angleB=158,linewidth=.4pt]{vl}{ryl}
   \nccurve[angleA=21,angleB=158,linewidth=.4pt]{wl}{rwl}
   \ncline[linewidth=.4pt]{xl}{rvl}

  \rput(2.05,2.8) {\scriptsize $w$} \rput(0.95,2.35) {\scriptsize $w_1$} \rput(3.15,2.3) {\scriptsize $w_2$}
   \rput(2.05,-1) {\scriptsize $ a$}  \rput(6.05,2.9) {\scriptsize $b$}
   \rput(2.05,-0.6) {$\textcolor[rgb]{1.00,0.00,0.00}{\times}$}  \rput(6.05,2.5) {$\textcolor[rgb]{1.00,0.00,0.00}{\times}$}

 \rput(4.05,-1.5){\scriptsize  $4D$ $THLN$($G_4$)}
\end{pspicture}
 \vskip50pt
\caption{
\label{f4}                                    
\footnotesize Example of weak vertex-pair}
\end{figure}

Unquestionably, for the weak vertex-pair $(w_1,w_2)$, any correct path $P_{w_1w_2}$ of length $l\geq3$ can't include the weak 2-degree vertex $w$.  It follows there is no correct hamiltonian path joining vertices $w_1$ and $w_2$ in $G_n-F$\cite{128}. However, we proved that $G_n-F(n\geq5)$ contains at most one weak $2$-degree vertex $w$ and one $w$-weak vertex-pair for any $F\subset V(G_n)\cup E(G_n)$ with $|F|\leq n-2$ in reference\cite{128}.
\begin{Def}$^{[15]}$
\textnormal{If $(w_1,w_2)$ is not a weak vertex-pair for any vertex $w\in V(G_n-F)$, then $(w_1,w_2)$ is addressed as a normal vertex pair.}
\end{Def}
\begin{thm}\label{lem9}$^{[15]}$\textnormal{Let  $F\subset V(G_n)\cup E(G_n)$ with $|F|\leq n-2$. Then for any vertex-pair $(u,v)$ in $G_n-F$, there is a $(n-2)$-fault-tolerant hamiltonian path $P_{uv}$ except $(u,v)$ being a weak vertex-pair.}
\end{thm}

\par In the paper, we studied the path-embedding in a $THLN$ with $n-2$ faulty elements and showed that if $F\subset V(G_n)\cup E(G_n)$ and $|F|\leq n-2$, then for arbitrary two different correct vertices $u$ and $v$, there is a fault-free path $P_{uv}$ of every length $l$ with $2^{n-1}-1\leq l\leq 2^n-f_v-1-\alpha$, where $\alpha=0$ if vertices $u$ and $v$ form a normal vertex-pair and $\alpha=1$ if vertices $u$ and $v$ form a weak vertex-pair in $G_n-F$($n\geq5$).

To do this simply, we can denote $G_n=L\oplus R$, where $L=G_{n-1}^0$ and $R=G_{n-1}^1$.  For any vertex $x\in L({\rm or}\ R)$, let $x^R({\rm or}\ x^L)$ be the sole vertex adjacent to vertex $x$ in $R({\rm or}\ L)$, and $N_L(x)( {\rm or}\ N_R(x))$ be the set of vertices that are adjacent to vertex $x$ in $L({\rm or}\ R)$. Let $E^C$ be the set of edges that join $L$ to $R$ and $E_L(x)({\rm or}\ E_R(x))$ be the set of edges incident to vertex $x$ in $L({\rm or}\ R)$.

 \par We use $P_{uv}$ to represent the path from vertex $u$ to vertex $v$. If $P_{uw}=(u,u_1,\cdots,u_s,w)$, $P_{wv}=(w,w_1,\cdots,w_t,v)$ and $V(P_{uw})\cap V(P_{wv})=\{w\}$, we use $P_{uw}+P_{wv}$ to denote the path $P_{uv}=(u,u_1,\cdots,u_s,w,w_1,\cdots,w_t,v)$, $P_{uv}(u_1,w_1)$ to represent the subpath of $P_{uv}$ which is from vertex $u_1$ to vertex $w_1$, $l_{uv}$ to denote the length of $P_{uv}$, $d_{uv}$ to denote the distance between vertex $u$ to vertex $v$. We denote $F^L=F\cap L$, $F^R=F\cap R$, $F^C=F\cap E^C$, $F_v=F\cap V(G_n)$, $F_e=F\cap E(G_n)$, $f_v=|F_v|$, $f_v^L=|F_v\cap V(L)|$, $f_v^R=|F_v\cap V(R)|$. We have $f_v=f_v^L+f_v^R$.

This paper is organized as below. Section 2 proved the main result. Section 3 concludes the paper.

\section{Main Result}
In the section, we will establish the main result of the paper. We depict theorem 2.1 as follows.

\begin{thm}\label{thm1}{If $F\subset V(G_n)\cup E(G_n)$ and $|F|\leq n-2$, then for any two distinct fault-free vertices $u$ and $v$, there exists a fault-free path $P_{uv}$ of every length $l$ with $2^{n-1}-1\leq l\leq 2^n-f_v-1-\alpha$, where $\alpha=0$ if vertices $u$ and $v$ form a normal vertex-pair and $\alpha=1$ if vertices $u$ and $v$ form a weak vertex-pair in $G_n-F$($n\geq5$).}
\end{thm}

\vskip 6pt

\noindent\textbf{Proof}.

 We prove the theorem by the induction on $n \geq 5$. The result holds
for $n=5$ by developing computer program using depth first searching
technique combining with backtracking and branch and bound algorithm. Assume that the theorem holds for $n-1$ with $n\geq6$, then we must show the theorem holds for $n$.
 In general, we assume $|F^R| \leq |F^L|$. Then $|F^R|\leq\lfloor \frac{n-2}{2}\rfloor\leq n-4$. Since for any vertex $x\in R$, $|N_R(x)|=n-1$. By $|F^R|\leq n-4$, $|N_{R-F^R}(x)|\geq3$. Then there is no weak vertex-pair in $R-F^R$.
 \par Let $u,v$ be any two distinct fault-free vertices in $G_n-F$. By Theorem 1.1, there is a faultless path $P_{uv}$ of length $l=2^n-f_v-1$ if vertices $u$ and $v$ form a normal vertex-pair in $G_n-F$. Then we only need to find each length $l$ with $2^{n-1}-1\leq l\leq2^n-f_v-2$ between arbitrary different vertices $u$ and $v$ in $G_n-F$. We divide the proof to two cases: (1), $|F^L| \leq n-3$; (2), $|F^L|= n-2$.
\par{\bf Case 1.} $|F^L| \leq n-3$.
\par{\bf Case 1.1.} $u,v\in V(L-F^L)$ or $u,v\in V(R-F^R)$. Firstly, We prove the case of $u,v\in V(L-F^L)$.
\par Since $|F^L| \leq n-3$, by induction hypothesis, there is a faultless path $P_{uv}$ of each length $l$ with $2^{n-2}-1\leq l\leq 2^{n-1}-f_v^L-2$ in $L-F^L$. Notice that there exist $\lfloor\frac{l+1}{2}\rfloor$ vertex-pairs in $P_{uv}$. Since $\lfloor\frac{l+1}{2}\rfloor-(n-2)\geq\frac{2^{n-2}}{2}-(n-2)\geq4(n\geq6)$, there is a faultless edge $ab\in E(P_{uv})$ with $a^R,b^R,aa^R,bb^R\notin F$. Since $|F^R|\leq n-4$, by induction hypothesis, there is a faultless path $P_{a^Rb^R}$ of each length $l_{a^Rb^R}$ with $2^{n-2}-1\leq l_{a^Rb^R}\leq 2^{n-1}-f_v^R-1$ in $R-F^R$. Let $P^1_{uv}=P_{uv}(u,a)+aa^R+P_{a^Rb^R}+b^Rb+P_{uv}(b,v)$. Then $P^1_{uv}$ is a faultless path of length $l^1_{uv}$ with $2^{n-1}-1\leq l^1_{uv}\leq 2^n-f_v-2$ in $G_n-F$(see Fig.3(a)).

For $u,v\in V(R-F^R)$, by a similar discussion, we can get a faultless path $P^1_{uv}$ of each length $l^1_{uv}$ with $2^{n-1}-1\leq l^1_{uv}\leq 2^n-f_v-2$ in $G_n-F$.

\par{\bf Case 1.2.} $u\in V(L-F^L)$ and $v\in V(R-F^R)$.
\par By the definition of $G_n$, $|E^C|=2^{n-1}$. Since $2^{n-1}-(n-2)\geq 28(n\geq6)$, there is a faultless edge $ab$ with $ab\in E^C$, $a,b\notin \{u,v\}$ and $a,b\notin F$.
 By induction hypothesis, there is a faultless path $P_{ua}$ of each length $l_{ua}$ with $2^{n-2}-1\leq l_{ua}\leq2^{n-1}-f_v^L-2$ in $L-F^L$ and a faultless path $P_{bv}$ of each length $l_{bv}$ with $2^{n-2}-1\leq l_{bv}\leq2^{n-1}-f_v^R-1$ in $R-F^R$. Let $P_{uv}=P_{ua}+ab+P_{bv}$. Then $P_{uv}$ is a faultless path of each length $l_{uv}$ with $2^{n-1}-1\leq l_{uv}\leq 2^n-f_v-2$ in $G_n-F$(see Fig.3(b)).
\begin{figure}[ht]
\psset{unit=0.9}

\begin{pspicture}(-3.65,-.3)(-1.15,2.7)

 \rput(.5,-.1){\rnode{a1}{}} \rput(.5,2){\rnode{b1}{}}
 \ncbox[nodesep=3pt,boxsize=1.,linearc=.3]{a1}{b1}
 \rput(3.5,-.1){\rnode{a3}{}} \rput(3.5,2){\rnode{b3}{}}
 \ncbox[nodesep=3pt,boxsize=1.,linearc=.3]{a3}{b3}
 \cnode(-.3,1.2){.09}{xl}\rput(-.1,1.2){\scriptsize$u$}
 \cnode(-.3,0.5){.09}{yl}\rput(-.1,.5){\scriptsize$v$}
 \cnode(1.3,1.2){.09}{al}\rput(1.1,1.2){\scriptsize$a$}
 \cnode(1.3,.5){.09}{bl}\rput(1.1,.5){\scriptsize$b$}
  \pnode(1.3,1.8){10}
  \pnode(-.3,1.8){21}
  \pnode(-.3,0.1){22}
  \ncline[linewidth=1.4pt]{xl}{21}
  \ncline[linewidth=1.4pt]{yl}{22}
  \ncline[linewidth=1.4pt]{al}{10}
 \pnode(1.3,0.1){14}
 \ncline[linewidth=1.4pt]{10}{21}
 \ncline[linewidth=1.4pt]{14}{22}
 \ncline[linewidth=1.4pt]{11}{12}
 \ncline[linewidth=1.4pt]{12}{14}
 \ncline[linewidth=1.4pt]{14}{bl}
 \ncline[linewidth=.5pt, linestyle=dashed]{al}{bl}

 \cnode(2.8,1.2){.09}{ar}\rput(3.15,1.25){\scriptsize$a^R$}
 \cnode(2.8,.5){.09}{br}\rput(3.15,.55){\scriptsize$b^R$}
 \ncline[linewidth=1.4pt]{al}{ar}
  \pnode(2.8,1.8){30}\ncline[linewidth=1.4pt]{ar}{30}
  \pnode(4.2,1.8){31}\ncline[linewidth=1.4pt]{30}{31}
  \pnode(4.2,0.1){32}\ncline[linewidth=1.4pt]{31}{32}
  \pnode(2.8,0.1){34}\ncline[linewidth=1.4pt]{32}{34}
  \ncline[linewidth=1.4pt]{34}{br}\ncline[linewidth=1.4pt]{bl}{br}
  \ncline[linewidth=.5pt, linestyle=dashed]{ar}{br}
 \rput(.5,2.5){\scriptsize$L=G^0_{n-1}$}\rput(3.5,2.5){\scriptsize$R=G^1_{n-1}$}
 \rput(2.,-.5){\scriptsize (a)}
 \end{pspicture}
\begin{pspicture}(-7.4,-.3)(-3.9,2.7)

 \rput(.5,-.1){\rnode{a1}{}} \rput(.5,2){\rnode{b1}{}}
 \ncbox[nodesep=3pt,boxsize=1.,linearc=.3]{a1}{b1}
 \rput(3.5,-.1){\rnode{a3}{}} \rput(3.5,2){\rnode{b3}{}}
 \ncbox[nodesep=3pt,boxsize=1.,linearc=.3]{a3}{b3}
 \pnode(-.3,1.8){30}
 \pnode(-.3,0.1){40}
 \ncline[linewidth=1.4pt]{30}{40}
 \cnode(1.3,1.3){.09}{u}\rput(1.05,1.3){\scriptsize$u$}
 \cnode(1.3,0.6){.09}{br}\rput(1.05,0.6){\scriptsize$a$}
  \pnode(1.3,1.8){10}
  \ncline[linewidth=1.4pt]{u}{10}
 \pnode(1.3,0.1){14}
  \ncline[linewidth=1.4pt]{30}{10}
 \ncline[linewidth=1.4pt]{40}{14}
 \ncline[linewidth=1.4pt]{14}{br}

 \cnode(2.8,1.3){.09}{ar}\rput(3.,1.3){\scriptsize$v$}
 \cnode(2.8,0.6){.09}{wr}\rput(3.,0.6){\scriptsize$b$}
 \ncline[linewidth=1.4pt]{br}{wr}
  \pnode(2.8,1.8){30}\ncline[linewidth=1.4pt]{ar}{30}
  \pnode(4.2,1.8){31}\ncline[linewidth=1.4pt]{30}{31}
  \pnode(4.2,0.1){32}\ncline[linewidth=1.4pt]{31}{32}
  \pnode(2.8,0.1){34}\ncline[linewidth=1.4pt]{32}{34}
  \ncline[linewidth=1.4pt]{34}{wr}
 \rput(.5,2.5){\scriptsize$L=G^0_{n-1}$}\rput(3.5,2.5){\scriptsize$R=G^1_{n-1}$}
 \rput(2.,-.5){\scriptsize (b)}
 \end{pspicture}
\caption{
\label{f4}                                    
\footnotesize  Illustrations of proofs of Case 1.1 and Case 1.2  of Theorem 2.1.}
\end{figure}
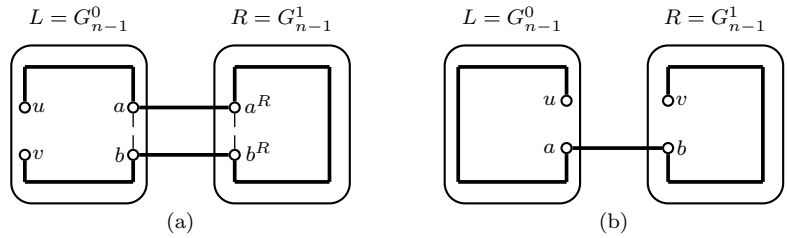
                          
\par{\bf Case 2.} $|F^L|= n-2$. Then $|F^R|=|F^C|=0$.
\par{\bf Case 2.1.} $|F^L\cap V(L)|\geq1$. Let $x\in F^L\cap V(L)$.
\par{\bf Case 2.1.1.} $u,v\in V(L-F^L)$.
\par We mark the faulty vertex $x$ as faultless temporarily. Let $F^L_1=F^L-x$, then $|F_1^L|=|F^L|-1=n-3$. By induction hypothesis, there is a faultless path $P_{uv}$ of each length $l_{uv}$ with $2^{n-2}-1\leq l_{uv}\leq 2^{n-1}-(f_v^L-1)-2=2^{n-1}-f_v^L-1$ in $L-F^L_{1}$. If the path $P_{uv}$ contains the faulty vertex $x$, let $a,b\in N_{P_{uv}}(x)$; otherwise, we can arbitrarily select a vertex $c$ from the path $P_{uv}$. Let $a,b\in N_{P_{uv}}(c)$. Since $|F^R|=0$, by induction hypothesis, there is a faultless path $P_{a^Rb^R}$ of each length $l_{a^Rb^R}$ with $2^{n-2}-1\leq l_{a^Rb^R}\leq 2^{n-1}-1$ in $R$. Let $P^1_{uv}=P_{uv}(u,a)+aa^R+P_{a^Rb^R}+b^Rb+P_{uv}(b,v)$. Then $P^1_{uv}$ is
a faultless path of each length $l^1_{uv}$ with $ 2^{n-1}-1\leq l^1_{uv}\leq 2^{n}-f_v-2$ in $G_n-F$(see Fig.4).
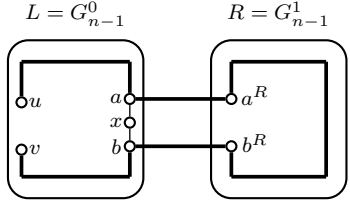
\begin{figure}[ht]
\psset{unit=0.9}
\begin{pspicture}(-6.6,-.3)(-3.1,2.7)
 \rput(.5,-.1){\rnode{a1}{}} \rput(.5,2){\rnode{b1}{}}
 \ncbox[nodesep=3pt,boxsize=1.,linearc=.3]{a1}{b1}
 \rput(3.5,-.1){\rnode{a3}{}} \rput(3.5,2){\rnode{b3}{}}
 \ncbox[nodesep=3pt,boxsize=1.,linearc=.3]{a3}{b3}

 \cnode(-.3,1.2){.09}{xl}\rput(-.1,1.2){\scriptsize$u$}
 \cnode(-.3,0.5){.09}{yl}\rput(-.1,.5){\scriptsize$v$}
 \cnode(1.3,1.25){.09}{al}\rput(1.1,1.25){\scriptsize$a$}
 \cnode(1.3,0.55){.09}{bl}\rput(1.1,.55){\scriptsize$b$}
 \cnode(1.3,0.9){.09}{bxl}\rput(1.1,.9){\scriptsize$x$}
  \pnode(1.3,1.8){10}
  \pnode(-.3,1.8){21}
  \pnode(-.3,0.1){22}
  \ncline[linewidth=1.4pt]{xl}{21}
  \ncline[linewidth=1.4pt]{yl}{22}
  \ncline[linewidth=1.4pt]{al}{10}
 \pnode(1.3,0.1){14}
 \ncline[linewidth=1.4pt]{10}{21}
 \ncline[linewidth=1.4pt]{14}{22}
 \ncline[linewidth=1.4pt]{11}{12}
 \ncline[linewidth=1.4pt]{12}{14}
 \ncline[linewidth=1.4pt]{14}{bl}
 \ncline[linewidth=.5pt, linestyle=dashed]{al}{bxl}
 \ncline[linewidth=.5pt, linestyle=dashed]{bl}{bxl}

 \cnode(2.8,1.25){.09}{ar}\rput(3.15,1.3){\scriptsize$a^R$}
 \cnode(2.8,0.55){.09}{br}\rput(3.15,0.6){\scriptsize$b^R$}
 \ncline[linewidth=1.4pt]{al}{ar}
  \pnode(2.8,1.8){30}\ncline[linewidth=1.4pt]{ar}{30}
  \pnode(4.2,1.8){31}\ncline[linewidth=1.4pt]{30}{31}
  \pnode(4.2,0.1){32}\ncline[linewidth=1.4pt]{31}{32}
  \pnode(2.8,0.1){34}\ncline[linewidth=1.4pt]{32}{34}
  \ncline[linewidth=1.4pt]{34}{br}\ncline[linewidth=1.4pt]{bl}{br}

 \rput(.5,2.5){\scriptsize$L=G^0_{n-1}$}\rput(3.5,2.5){\scriptsize$R=G^1_{n-1}$}

 \end{pspicture}

\caption{
\label{f4}                                  
\footnotesize  Illustrations of proofs of Case 2.1.1  of Theorem 2.1.}
\end{figure}                           

\par{\bf Case 2.1.2.} $u\in V(L-F^L)$ and $v\in V(R)$.
\par  We mark the faulty vertex $x$ as faultless temporarily. Let $F^L_1=F^L-x$, then $|F^L_1|=|F^L|-1=n-3$. By induction hypothesis, there is a faultless path $P_{ux}$ of each length $l_{ux}$ with $2^{n-2}-1\leq l_{ux}\leq2^{n-1}-(f_v^L-1)-2=2^{n-1}-f_v^L-1$ in $L-F^L_1$. Let $x_1\in N_{P_{ux}}(x)$.
\par{\bf Case 2.1.2.1.} $x_1^R=v$.
 \par Let $ab\in E(P_{ux})$ with $a,b\notin\{u,x_1,x\}$. We mark the correct vertex $v$ as faulty temporarily. Let $F^R_1=F^R+v$, then $|F^R_1|=|F^R|+1\leq n-4(n\geq6)$. By induction hypothesis, there is a faultless path $P_{a^Rb^R}$ of each length $l_{a^Rb^R}$ with $2^{n-2}-1\leq l_{a^Rb^R}\leq2^{n-1}-2$ in $R-F^R_1$. Let $P_{uv}=P_{ux}(u,a)+aa^R+P_{a^Rb^R}+b^Rb+P_{ux}(b,x_1)+x_1v$. Then $P_{uv}$ is a faultless path of each length $l_{uv}$ with $2^{n-1}-1\leq l_{uv}\leq 2^n-f_v-2$ in $G_n-F$(see Fig.5(a)).
\par{\bf Case 2.1.2.2.} $x_1^R\neq v$.
 \par By induction hypothesis, there is a faultless path $P_{x_1^Rv}$ of each length $l_{x_1^Rv}$ with $2^{n-2}-1\leq l_{x_1^Rv}\leq 2^{n-1}-1$ in $R$. Let $P_{uv}=P_{ux}(u,x_1)+x_1x_1^R+P_{x_1^Rv}$. Then $P_{uv}$ is a faultless path of each length $l_{uv}$ with $2^{n-1}-1\leq l_{uv}\leq 2^n-f_v-2$ in $G_n-F$(see Fig.5(b)).
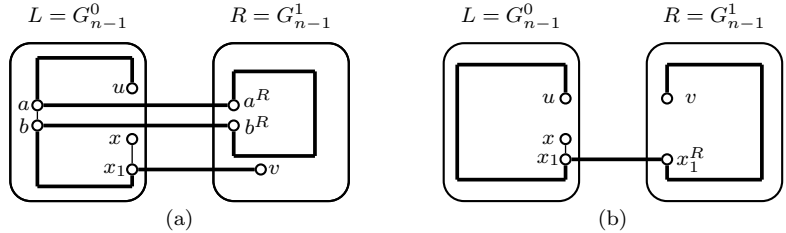
\begin{figure}[ht]
\psset{unit=0.9}
\begin{pspicture}(-3.65,-.3)(-1.15,2.7)
 \rput(.5,-.1){\rnode{a1}{}} \rput(.5,2){\rnode{b1}{}}
 \ncbox[nodesep=3pt,boxsize=1.,linearc=.3]{a1}{b1}
 \rput(3.5,-.1){\rnode{a3}{}} \rput(3.5,2){\rnode{b3}{}}
 \ncbox[nodesep=3pt,boxsize=1.,linearc=.3]{a3}{b3}
\rput(.5,-.1){\rnode{a1}{}} \rput(.5,2){\rnode{b1}{}}
 \ncbox[nodesep=3pt,boxsize=1.,linearc=.3]{a1}{b1}
 \rput(3.5,-.1){\rnode{a3}{}} \rput(3.5,2){\rnode{b3}{}}
 \ncbox[nodesep=3pt,boxsize=1.,linearc=.3]{a3}{b3}
 \cnode(1.3,1.45){.09}{al}\rput(1.1,1.45){\scriptsize$u$}
 \cnode(1.3,0.25){.09}{bl}\rput(1.05,.25){\scriptsize$x_1$}
 \cnode(3.2,0.25){.09}{blvv}\rput(3.4,.25){\scriptsize$v$}
  \ncline[linewidth=1.4pt]{bl}{blvv}
 \cnode(-.1,1.2){.09}{xl}\rput(-.3,1.2){\scriptsize$a$}
 \cnode(-.1,0.9){.09}{yl}\rput(-.3,.9){\scriptsize$b$}
 \cnode(1.3,0.7){.09}{yxl}\rput(1.05,.7){\scriptsize$x$}
  \pnode(1.3,1.9){10}
  \pnode(-.1,1.9){21}
  \pnode(-.1,-0.0){22}
  \ncline[linewidth=1.4pt]{xl}{21}
  \ncline[linewidth=1.4pt]{yl}{22}
  \ncline[linewidth=1.4pt]{al}{10}
 \pnode(1.3,-0.0){14}
 \ncline[linewidth=1.4pt]{10}{21}
 \ncline[linewidth=1.4pt]{14}{22}
 \ncline[linewidth=1.4pt]{11}{12}
 \ncline[linewidth=1.4pt]{12}{14}
 \ncline[linewidth=1.4pt]{14}{bl}
 \ncline[linewidth=.5pt, linestyle=dashed]{xl}{yl}
 \ncline[linewidth=.5pt, linestyle=dashed]{bl}{yxl}

 \cnode(2.8,1.2){.09}{ar}\rput(3.15,1.3){\scriptsize$a^R$}
 \cnode(2.8,0.9){.09}{br}\rput(3.15,0.9){\scriptsize$b^R$}
 \ncline[linewidth=1.4pt]{xl}{ar}
  \pnode(2.8,1.7){30}\ncline[linewidth=1.4pt]{ar}{30}
  \pnode(4.,1.7){31}\ncline[linewidth=1.4pt]{30}{31}
  \pnode(4.,0.45){32}\ncline[linewidth=1.4pt]{31}{32}
  \pnode(2.8,0.45){34}\ncline[linewidth=1.4pt]{32}{34}
  \ncline[linewidth=1.4pt]{34}{br}\ncline[linewidth=1.4pt]{yl}{br}
 \rput(.5,2.5){\scriptsize$L=G^0_{n-1}$}\rput(3.5,2.5){\scriptsize$R=G^1_{n-1}$}
 \rput(2.,-.5){\scriptsize (a)}
 \end{pspicture}
\begin{pspicture}(-7.4,-.3)(-3.9,2.7)

\rput(.5,-.1){\rnode{a1}{}} \rput(.5,2){\rnode{b1}{}}
 \ncbox[nodesep=3pt,boxsize=1.,linearc=.3]{a1}{b1}
 \rput(3.5,-.1){\rnode{a3}{}} \rput(3.5,2){\rnode{b3}{}}
 \ncbox[nodesep=3pt,boxsize=1.,linearc=.3]{a3}{b3}
 \pnode(-.3,1.8){30}
 \pnode(-.3,0.1){40}
 \ncline[linewidth=1.4pt]{30}{40}
 \cnode(1.3,1.3){.09}{u}\rput(1.05,1.3){\scriptsize$u$}
 \cnode(1.3,0.7){.09}{bl}\rput(1.05,0.7){\scriptsize$x$}
 \cnode(1.3,0.4){.09}{br}\rput(1.05,0.4){\scriptsize$x_1$}
  \ncline[linewidth=.5pt, linestyle=dashed]{bl}{br}
  \pnode(1.3,1.8){10}
  \ncline[linewidth=1.4pt]{u}{10}
 \pnode(1.3,0.1){14}
  \ncline[linewidth=1.4pt]{30}{10}
 \ncline[linewidth=1.4pt]{40}{14}
 \ncline[linewidth=1.4pt]{14}{br}

 \cnode(2.8,1.3){.09}{arrr}\rput(3.15,1.3){\scriptsize$v$}
 \cnode(2.8,0.4){.09}{wr}\rput(3.15,0.4){\scriptsize$x_1^R$}
 \ncline[linewidth=1.4pt]{br}{wr}
  \pnode(2.8,1.8){30}\ncline[linewidth=1.4pt]{arrr}{30}
  \pnode(4.2,1.8){31}\ncline[linewidth=1.4pt]{30}{31}
  \pnode(4.2,0.1){32}\ncline[linewidth=1.4pt]{31}{32}
  \pnode(2.8,0.1){34}\ncline[linewidth=1.4pt]{32}{34}
  \ncline[linewidth=1.4pt]{34}{wr}
  \rput(.5,2.5){\scriptsize$L=G^0_{n-1}$}\rput(3.5,2.5){\scriptsize$R=G^1_{n-1}$}
 \rput(2.,-.5){\scriptsize (b)}
 \end{pspicture}

\caption{
\label{f4}                                    
\footnotesize  Illustrations of proofs of  Case 2.1.2  of Theorem 2.1.}
\end{figure}                                   

\par{\bf Case 2.1.3.} $u,v\in V(R)$.
\par Since $|F^R|=0$, by induction hypothesis, there is a faultless path $P_{uv}$ of length $l=2^{n-1}-1$ in $R$. Thus, we only need to consider each length $l$ with  $2^{n-1}\leq l_{uv}\leq 2^{n}-f_v-2$.
\par{\bf Case 2.1.3.1.} $|\{u^L,v^L\}\cap F_v|\geq1$. In general, assume $u^L\in F_v$. We mark the faulty vertex $u^L$ as faultless temporarily. Let $F^L_1=F^L-u^L$, then $|F^L_1|=|F^L|-1=n-3$.
\par Let $S=N_R(v)-u$. Then $|S|\geq n-2$. Since $|F^L_1|=n-3$, there is a vertex $v_1\in S$ with $v_1^L\notin F$. By induction hypothesis, there is a faultless path $P_{u^Lv_1^L}$ of each length $l_{u^Lv_1^L}$ with $2^{n-2}-1\leq l_{u^Lv_1^L}\leq 2^{n-1}-(f_v^L-1)-2=2^{n-1}-f_v^L-1$ in $L-F_1^L$. Let $u_1\in N_{P_{u^Lv_1^L}}(u^L)$.
\par If $u_1^R\neq v$, let $F_1^R=F^R+\{v_1,v\}$, then $|F_1^R|=|F^R|+2=2\leq n-4(n\geq6)$. By induction hypothesis, there is a faultless path $P_{uu_1^R}$ of each length $l_{uu_1^R}$ with $2^{n-2}-1\leq l_{uu_1^R}\leq 2^{n-1}-3$ in $R-F_1^R$. Let $P_{uv}=P_{uu_1^R}+u_1^Ru_1+P_{u^Lv_1^L}(u_1,v_1^L)+v_1^Lv_1+v_1v$. Then $P_{uv}$ is a faultless path of each length $l_{uv}$ with $2^{n-1}\leq l_{uv}\leq 2^{n}-f_v-2$ in $G_n-F$(See Fig.6(a)).
\par If $u_1^R=v$, let $F^R_1=F^R+v$, then $|F^R_1|=1\leq n-4(n\geq6)$. By induction hypothesis, there is a faultless path $P_{uv_1}$ of each length $l_{uv_1}$ with $2^{n-2}-1\leq l_{uv_1}\leq 2^{n-1}-2$ in $R-F^R_1$. Let $P_{uv}=P_{uv_1}+v_1v^L_1+P_{u^Lv_1^L}(v_1^L, u_1)+u_1v$. Then $P_{uv}$ is a fault-free path of each length $l_{uv}$ with $2^{n-1} \leq l_{uv}\leq 2^{n}-f_v-2$ in $G_n-F$(See Fig.6(b)).
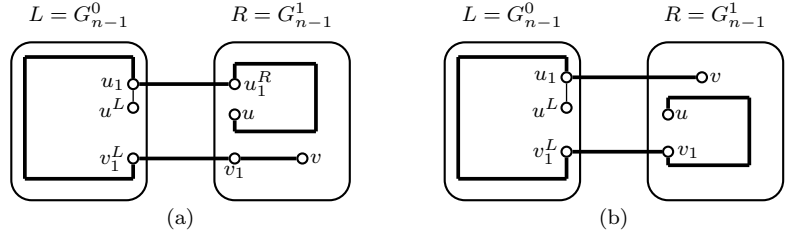
\begin{figure}[ht]
\psset{unit=0.9}

\begin{pspicture}(-3.65,-.3)(-1.15,2.7)
 \rput(.5,-.1){\rnode{a1}{}} \rput(.5,2){\rnode{b1}{}}
 \ncbox[nodesep=3pt,boxsize=1.,linearc=.3]{a1}{b1}
 \rput(3.5,-.1){\rnode{a3}{}} \rput(3.5,2){\rnode{b3}{}}
 \ncbox[nodesep=3pt,boxsize=1.,linearc=.3]{a3}{b3}

  \cnode(1.3,1.5){.09}{al}\rput(1.,1.5){\scriptsize$u_1$}
 \cnode(1.3,1.15){.09}{u}\rput(1.,1.15){\scriptsize$u^L$}
 \ncline[linewidth=.5pt, linestyle=dashed]{al}{u}
 \cnode(1.3,0.4){.09}{brx}\rput(1.,0.4){\scriptsize$v_1^L$}
 \cnode(3.8,0.4){.09}{av}\rput(4.0,0.4){\scriptsize$v$}
 \cnode(2.8,1.5){.09}{alm}\rput(3.15,1.5){\scriptsize$u_1^R$}
 \cnode(2.8,1.05){.09}{un}\rput(3.0,1.05){\scriptsize$u$}
 \cnode(2.8,0.4){.09}{blb}\rput(2.8,0.2){\scriptsize$v_1$}
  \ncline[linewidth=1.4pt]{brx}{blb}
  \ncline[linewidth=1.4pt]{av}{blb}
    \ncline[linewidth=1.4pt]{al}{alm}

  \pnode(1.3,1.9){10}
  \pnode(-.3,1.9){21}
  \pnode(-.3,0.1){22}
  \ncline[linewidth=1.4pt]{al}{10}
  \ncline[linewidth=1.4pt]{21}{22}
 \pnode(1.3,0.1){14}
 \ncline[linewidth=1.4pt]{10}{21}
 \ncline[linewidth=1.4pt]{14}{22}
 \ncline[linewidth=1.4pt]{11}{12}
 \ncline[linewidth=1.4pt]{12}{14}
 \ncline[linewidth=1.4pt]{14}{brx}
 \ncline[linewidth=.5pt, linestyle=dashed]{xl}{yl}
 \ncline[linewidth=1.4pt]{xl}{ar}
  \pnode(2.8,1.8){30} \ncline[linewidth=1.4pt]{alm}{30}
  \pnode(4.,1.8){31}\ncline[linewidth=1.4pt]{30}{31}
  \pnode(4.,0.8){32}\ncline[linewidth=1.4pt]{31}{32}
  \pnode(2.8,0.8){34}\ncline[linewidth=1.4pt]{32}{34}

 \ncline[linewidth=1.4pt]{un}{34}

 \rput(.5,2.5){\scriptsize$L=G^0_{n-1}$}\rput(3.5,2.5){\scriptsize$R=G^1_{n-1}$}
 \rput(2.,-.5){\scriptsize (a)}
 \end{pspicture}
\begin{pspicture}(-7.4,-.3)(-3.9,2.7)
 \rput(.5,-.1){\rnode{a1}{}} \rput(.5,2){\rnode{b1}{}}
 \ncbox[nodesep=3pt,boxsize=1.,linearc=.3]{a1}{b1}
 \rput(3.5,-.1){\rnode{a3}{}} \rput(3.5,2){\rnode{b3}{}}
 \ncbox[nodesep=3pt,boxsize=1.,linearc=.3]{a3}{b3}

 \cnode(1.3,1.6){.09}{al}\rput(1.,1.6){\scriptsize$u_1$}
 \cnode(1.3,1.15){.09}{u}\rput(1.,1.15){\scriptsize$u^L$}
 \ncline[linewidth=.5pt, linestyle=dashed]{al}{u}
 \cnode(1.3,0.5){.09}{brx}\rput(1.,0.5){\scriptsize$v_1^L$}
 \cnode(3.3,1.6){.09}{alm}\rput(3.5,1.6){\scriptsize$v$}
 \cnode(2.8,1.05){.09}{un}\rput(3.0,1.05){\scriptsize$u$}
 \cnode(2.8,0.5){.09}{blb}\rput(3.1,0.5){\scriptsize$v_1$}
  \ncline[linewidth=1.4pt]{brx}{blb}
  \ncline[linewidth=1.4pt]{al}{alm}

  \pnode(1.3,1.9){10}
  \pnode(-.3,1.9){21}
  \pnode(-.3,0.1){22}
  \ncline[linewidth=1.4pt]{al}{10}
  \ncline[linewidth=1.4pt]{21}{22}
 \pnode(1.3,0.1){14}
 \ncline[linewidth=1.4pt]{10}{21}
 \ncline[linewidth=1.4pt]{14}{22}
 \ncline[linewidth=1.4pt]{11}{12}
 \ncline[linewidth=1.4pt]{12}{14}
 \ncline[linewidth=1.4pt]{14}{brx}
 \ncline[linewidth=.5pt, linestyle=dashed]{xl}{yl}
 \ncline[linewidth=1.4pt]{xl}{ar}
  \pnode(2.8,1.3){30} \ncline[linewidth=1.4pt]{un}{30}
  \pnode(4.,1.3){31}\ncline[linewidth=1.4pt]{30}{31}
  \pnode(4.,0.3){32}\ncline[linewidth=1.4pt]{31}{32}
  \pnode(2.8,0.3){34}\ncline[linewidth=1.4pt]{32}{34}

 \ncline[linewidth=1.4pt]{blb}{34}
 \rput(.5,2.5){\scriptsize$L=G^0_{n-1}$}\rput(3.5,2.5){\scriptsize$R=G^1_{n-1}$}
 \rput(2.,-.5){\scriptsize (b)}
 \end{pspicture}
\caption{
\label{f4}                                   
\footnotesize  Illustrations of proofs of Case 2.1.3.1.  of Theorem 2.1.}
\end{figure}                         

\par{\bf Case 2.1.3.2.} $|\{u^L,v^L\}\cap F_v|=0$.  We mark the faulty vertex $x$ as faultless temporarily. Let $F^L_1=F^L-x$, then $|F^L_1|=|F^L|-1=n-3$.

\par Since $|F^L_1|=n-3$, by induction hypothesis, there is a faultless path $P_{xv^L}$ of each length $2^{n-2}-1\leq l_{xv^L}\leq 2^{n-1}-(f_v^L-1)-2=2^{n-1}-f_v^L-1$ in $L-F_1^L$. Let $x_1\in N_{P_{xv^L}}(x)$.
\par If $x_1^R=u$, let $F_1^R=F^R+\{u,v\}$, then $|F_1^R|=|F^R|+2=2\leq n-4(n\geq6)$. Let $ab\in E(P_{xv^L})$ with $a,b\notin\{x,x_1,v^L\}$. By induction hypothesis, there is a faultless path $P_{a^Rb^R}$ of each length $l_{a^Rb^R}$ with $2^{n-2}-1 \leq l_{a^Rb^R}\leq 2^{n-1}-3$ in $R-F_1^R$.  Let $P_{uv}=ux_1+P_{xv^L}(x_1,a)+aa^R+P_{a^Rb^R}+b^Rb+P_{xv^L}(b,v^L)+v^Lv$. Then $P_{uv}$ is a faultless path of each length $l_{uv}$ with $2^{n-1}\leq l_{uv}\leq 2^n-f_v-2$ in $G_n-F$(See Fig.7(a)).
\par If $x_1^R\neq u$, let $F_1^R=F^R+v$, then $|F_1^R|=|F^R|+1=1\leq n-4(n\geq6)$. By induction hypothesis, there is a faultless path $P_{ux_1^R}$ of each length $l_{ux_1^R}$ with $2^{n-2}-1\leq l_{ux_1^R}\leq 2^{n-1}-2$ in $R-F_1^R$. Let $P_{uv}=P_{ux_1^R}+x_1^Rx_1+P_{xv^L}(x_1,v^L)+v^Lv$. Then $P_{uv}$ is a
faultless path of each length $l_{uv}$ with $2^{n-1}\leq l_{uv}\leq 2^n-f_v-2$ in $G_n-F$(See Fig.7(b)).
 \begin{figure}[ht]
\psset{unit=0.9}
\begin{pspicture}(-3.65,-.3)(-1.15,2.7)
 \rput(.5,-.1){\rnode{a1}{}} \rput(.5,2){\rnode{b1}{}}
 \ncbox[nodesep=3pt,boxsize=1.,linearc=.3]{a1}{b1}
 \rput(3.5,-.1){\rnode{a3}{}} \rput(3.5,2){\rnode{b3}{}}
 \ncbox[nodesep=3pt,boxsize=1.,linearc=.3]{a3}{b3}

 \cnode(1.3,1.6){.09}{al}\rput(1.,1.6){\scriptsize$x_1$}
  \cnode(1.3,1.25){.09}{alx}\rput(1.,1.25){\scriptsize$x$}
 \cnode(1.3,0.2){.09}{bl}\rput(1.,0.2){\scriptsize$v^L$}
 \cnode(3.5,0.2){.09}{blvvv}\rput(3.7,0.2){\scriptsize$v$}
 \cnode(3.5,1.6){.09}{blvvu}\rput(3.7,1.6){\scriptsize$u$}
 \ncline[linewidth=1.4pt]{al}{blvvu}
  \ncline[linewidth=1.4pt]{bl}{blvvv}
 \cnode(-.1,1.0){.09}{xl}\rput(-.32,1.0){\scriptsize$a$}
 \cnode(-.1,0.6){.09}{yl}\rput(-.32,.6){\scriptsize$b$}
  \pnode(1.3,1.9){10}
  \pnode(-.1,1.9){21}
  \pnode(-.1,0.){22}
  \ncline[linewidth=1.4pt]{xl}{21}
  \ncline[linewidth=1.4pt]{yl}{22}
  \ncline[linewidth=1.4pt]{al}{10}
  \ncline[linewidth=.5pt, linestyle=dashed]{al}{alx}
 \pnode(1.3,0.){14}
 \ncline[linewidth=1.4pt]{10}{21}
 \ncline[linewidth=1.4pt]{14}{22}
 \ncline[linewidth=1.4pt]{11}{12}
 \ncline[linewidth=1.4pt]{12}{14}
 \ncline[linewidth=1.4pt]{14}{bl}
 \ncline[linewidth=.5pt, linestyle=dashed]{xl}{yl}
 
 \cnode(2.8,1.0){.09}{arr}\rput(3.15,1.1){\scriptsize$a^R$}
 \cnode(2.8,0.6){.09}{br}\rput(3.15,0.6){\scriptsize$b^R$}
 \ncline[linewidth=1.4pt]{xl}{arr}
  \pnode(2.8,1.35){30}\ncline[linewidth=1.4pt]{arr}{30}
  \pnode(3.8,1.35){31}\ncline[linewidth=1.4pt]{30}{31}
  \pnode(3.8,0.4){32}\ncline[linewidth=1.4pt]{31}{32}
  \pnode(2.8,0.4){34}\ncline[linewidth=1.4pt]{32}{34}
  \ncline[linewidth=1.4pt]{34}{br}\ncline[linewidth=1.4pt]{yl}{br}
  \ncline[linewidth=.5pt, linestyle=dashed]{ar}{br}

 \rput(.5,2.5){\scriptsize$L=G^0_{n-1}$}\rput(3.5,2.5){\scriptsize$R=G^1_{n-1}$}
 \rput(2.,-.5){\scriptsize (a)}
 \end{pspicture}
\begin{pspicture}(-7.4,-.3)(-3.9,2.7)
 \rput(.5,-.1){\rnode{a1}{}} \rput(.5,2){\rnode{b1}{}}
 \ncbox[nodesep=3pt,boxsize=1.,linearc=.3]{a1}{b1}
 \rput(3.5,-.1){\rnode{a3}{}} \rput(3.5,2){\rnode{b3}{}}
 \ncbox[nodesep=3pt,boxsize=1.,linearc=.3]{a3}{b3}

 \cnode(1.3,1.5){.09}{al}\rput(1.,1.5){\scriptsize$x_1$}
 \cnode(1.3,1.15){.09}{u}\rput(1.,1.15){\scriptsize$x$}
 \ncline[linewidth=.5pt, linestyle=dashed]{al}{u}
 \cnode(1.3,0.5){.09}{brx}\rput(1.,0.5){\scriptsize$v^L$}
 \cnode(2.8,1.5){.09}{alm}\rput(3.15,1.5){\scriptsize$x_1^R$}
 \cnode(2.8,1.05){.09}{un}\rput(3.05,1.05){\scriptsize$u$}
 \cnode(3.1,0.5){.09}{blb}\rput(3.3,0.5){\scriptsize$v$}
  \ncline[linewidth=1.4pt]{brx}{blb}
    \ncline[linewidth=1.4pt]{al}{alm}

  \pnode(1.3,1.9){10}
  \pnode(-.3,1.9){21}
  \pnode(-.3,0.1){22}
  \ncline[linewidth=1.4pt]{al}{10}
  \ncline[linewidth=1.4pt]{21}{22}
 \pnode(1.3,0.1){14}
 \ncline[linewidth=1.4pt]{10}{21}
 \ncline[linewidth=1.4pt]{14}{22}
 \ncline[linewidth=1.4pt]{11}{12}

 \ncline[linewidth=1.4pt]{12}{14}
 \ncline[linewidth=1.4pt]{14}{brx}
 \ncline[linewidth=.5pt, linestyle=dashed]{xl}{yl}
 \ncline[linewidth=1.4pt]{xl}{ar}
  \pnode(2.8,1.8){30} \ncline[linewidth=1.4pt]{alm}{30}
  \pnode(4.,1.8){31}\ncline[linewidth=1.4pt]{30}{31}
  \pnode(4.,0.8){32}\ncline[linewidth=1.4pt]{31}{32}
  \pnode(2.8,0.8){34}\ncline[linewidth=1.4pt]{32}{34}
   \ncline[linewidth=1.4pt]{un}{34}
 \rput(.5,2.5){\scriptsize$L=G^0_{n-1}$}\rput(3.5,2.5){\scriptsize$R=G^1_{n-1}$}
 \rput(2.,-.5){\scriptsize (b)}
 \end{pspicture}
\caption{
\label{f4}                                    
\footnotesize  Illustrations of proofs of Case 2.1.3.2.  of Theorem 2.1.}
\end{figure}
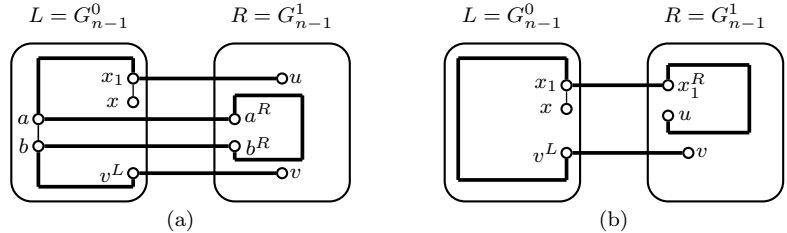                        

\par{\bf Case 2.2.} $|F^L\cap V(L)|=0$. Then $F^L=F_e=F$.
\par Let $x\in V(G_n)$ with $d_{G_n-F}(x)=\delta(G_n-F)$.
\par{\bf Case 2.2.1.} $x\notin \{u,v\}$.
\par Let $e$ be an edge with $e\in F_e$,  $F^1=F-e+\{x\}$, then $|F^1|=n-2$ and $|F^1_v|=1$. We show that $(u,v)$ is a normal vertex pair in $G_n-F^1$ as follows.
\par If $\delta(G_n-F)\geq 4$, we discuss $\delta(G_n-F^1)$ in the following four cases.

\par \indent(1) For any correct vertex $x_1\in N_{G_n-F}(x)$ with $e\notin E_{G_n}(x_1)$. Notice that $\delta(G_n-F)\geq 4$, then $d_{G_n-F^1}(x_1)=d_{G_n-F}(x_1)-1\geq3$.
\par \indent(2) For any correct vertex $x_1\in N_{G_n-F}(x)$ with $e\in E_{G_n}(x_1)$. Since $F^1=F-e+\{x\}$, we have $d_{G_n-F^1}(x_1)=d_{G_n-F}(x_1)\geq4$.
\par \indent(3) For any correct vertex $x_1\notin N_{G_n-F}(x)$ with $e\in E_{G_n}(x_1)$. Notice that $\delta(G_n-F)\geq 4$, then $d_{G_n-F^1}(x_1)=d_{G_n-F}(x_1)+1\geq5$.
\par \indent(4) For any correct vertex $x_1\notin N_{G_n-F}(x)$ with $e\notin E_{G_n}(x_1)$, Since $F^1=F-e+\{x\}$, we have $d_{G_n-F^1}(x_1)=d_{G_n-F}(x_1)\geq4$.
\par Above all, we conclude that $\delta(G_n-F^1)\geq 3$.
\par If $\delta(G_n-F)\leq 3$, then $|E_{G_n}(x)\cap F|\geq n-3$. For any $z\in V(G_n-F^1)$, since $|F|=n-2$ and $|E_{G_n}(x)\cap E_{G_n}(z)|\leq1$, we have $|(E_{G_n}(z)\cup N_{G_n}(z))\cap F^1|\leq 3$. It follows that $\delta(G_n-F^1)\geq n-3\geq3(n\geq6)$.
\par Hence, there is no weak vertex-pair in $G_n-F^1$, \emph{i.e.}, $(u,v)$ is a normal  vertex pair in $G_n-F^1$. By the proof of Case 2.1 and Theorem 1.1, there is a faultless path $P_{uv}$ of every length $l$ with $2^{n-1}-1\leq l\leq 2^n-|F^1_v|-1=2^n-2$ in $G_n-F$($n\geq5$).
\par{\bf Case 2.2.2.} $x\in \{u,v\}$. In general, assume that $x=u$.

\par Let $e$ be an edge with $e=uy\in F_e$ with $y\neq v$ and $F^1=F-e+\{y\}$, then $|F^1|=n-2$ and $|F^1_v|=1$. We show that $(u,v)$ is a normal vertex pair in $G_n-F^1$ as follows.
\par Let $z$ be an arbitrary vertex of $V(G_n-F^1)-\{u,v\}$.
\par If $\delta(G_n-F)\geq 4$, similar to the above discussion in Case 2.2.1, we have $\delta(G_n-F^1)\geq 3$. It means that $d_{G_n-F^1}(z)\geq3$.
\par If $\delta(G_n-F)\leq 3$, then $|E_{G_n}(u)\cap F|\geq n-3$. Since $|F|=n-2$ and $|E_{G_n}(u)\cap E_{G_n}(z)|\leq1$, we have $|(E_{G_n}(z)\cup N_{G_n}(z))\cap F^1|\leq 3$. It follows that $d_{G_n-F^1}(z)\geq n-3\geq3(n\geq6)$.
\par Hence, $(u,v)$ can not be a $z$-weak vertex pair in $G_n-F^1$, i.e., $(u,v)$ is a normal vertex pair in $G_n-F^1$. By the proof of Case 2.1 and Theorem 1.1, there is a faultless path $P_{uv}$ of every length $l$ with $2^{n-1}-1\leq l\leq 2^n-|F^1_v|-1=2^n-|F^1_v|-1=2^n-2$ in $G_n-F$($n\geq5$).
\qed

\section{Concluding Remarks}
This paper considered the path-embedding in an $n$-$D$ $THLN$($n\geq5$) with a set $F$ of up to $n-2$ faulty elements. We have proved that for arbitrary two different correct vertices $u$ and $v$, there exists a fault-free path $P_{uv}$ of every length $l$ with $2^{n-1}-1\leq l\leq 2^n-f_v-1-\alpha$, where $\alpha=0$ if vertices $u$ and $v$ form a normal vertex-pair and $\alpha=1$ if vertices $u$ and $v$ form a weak vertex-pair in $G_n-F$($n\geq5$).
The proposed theorem in the paper can be applied to several multiprocessor systems, including $n$-dimensional M\"{o}bius cubes $MQ_n$\cite{2}, $n$-dimensional locally twisted cubes $LTQ_n$\cite{24}, $n$-dimensional twisted cubes $TQ_n$\cite{15} for odd $n$, and $n$-dimensional crossed cubes $CQ_n$\cite{6}. Fig.8 illustrates $CQ_3$ and $CQ_4$. The graphs shown in Fig.9 are $LTQ_3$ and $LTQ_4$. Fig.10 plots $MQ^0_3$, $MQ^0_4$, $MQ^1_3$ and $MQ^1_4$. The graphs shown in Fig.11 are $TQ_3$ and $TQ_5$. By the discussion in reference\cite{128}, $MQ_n,LTQ_n,TQ_n,CQ_n\in THLNs$.
\par Hence, by Theorem 2.1, we can obtain the following four Corollaries.
\begin{cor}\textnormal{If $F\subset V(CQ_n)\cup E(CQ_n)$ with $|F|\leq n-2$, then for any two different correct vertices $u$ and $v$, there is a fault-free path $P_{uv}$ of every length $l$ with $2^{n-1}-1\leq l\leq 2^n-f_v-1-\alpha$, where $\alpha=0$ if vertices $u$ and $v$ form a normal vertex-pair and $\alpha=1$ if vertices $u$ and $v$ form a weak vertex-pair in $CQ_n-F$($n\geq5$).}
\end{cor}
\begin{cor}\textnormal{If $F\subset V(LTQ_n)\cup E(LTQ_n)$ with $|F|\leq n-2$, then for any two different correct vertices $u$ and $v$, there exists a fault-free path $P_{uv}$ of every length $l$ with $2^{n-1}-1\leq l\leq 2^n-f_v-1-\alpha$, where $\alpha=0$ if vertices $u$ and $v$ form a normal vertex-pair and $\alpha=1$ if vertices $u$ and $v$ form a weak vertex-pair in $LTQ_n-F$($n\geq5$). }
\end{cor}
\vskip10pt
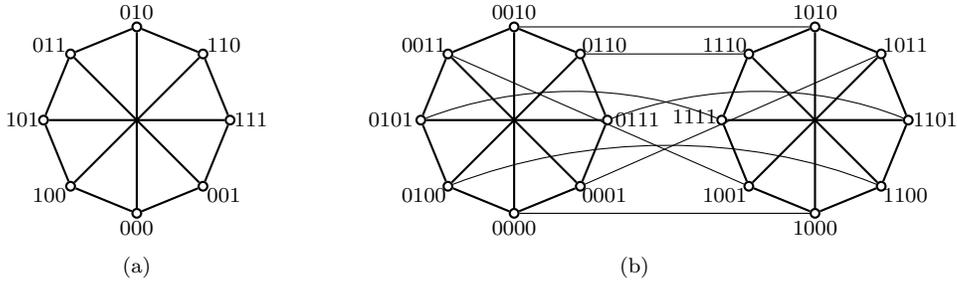
\begin{figure}[ht]
\psset{unit=0.8}
\begin{pspicture}(-2,1)(0.5,2.2)

  \cnode(0.95,2.05){.09}{xl}\rput(0.6,2.2){\scriptsize$011$}
  \cnode(0.95,-.15){.09}{vl}\rput(0.6,-.3){\scriptsize$100$}
  \cnode(3.15,2.05){.09}{ul}\rput(3.5,2.2){\scriptsize$110$}
  \cnode(3.15,-.15){.09}{yl}\rput(3.5,-0.3){\scriptsize$001$}
  \cnode(0.5,0.95){.09}{wl} \rput(0.15,0.95){\scriptsize$101$}
  \cnode(3.6,0.95){.09}{zl} \rput(3.95,0.95){\scriptsize$111$}
  \cnode(2.05,2.5){.09}{wla}\rput(2.05,2.75){\scriptsize$010$}
  \cnode(2.05,-0.6){.09}{zlb}\rput(2.05,-0.85){\scriptsize$000$}

  \ncline[linewidth=0.8pt]{xl}{wla}
  \ncline[linewidth=0.8pt]{wl}{xl}
   \ncline[linewidth=0.8pt]{yl}{xl}
   \ncline[linewidth=0.8pt]{zlb}{wla}
   \ncline[linewidth=0.8pt]{ul}{wla}
   \ncline[linewidth=0.8pt]{ul}{zl}
   \ncline[linewidth=0.8pt]{ul}{vl}
   \ncline[linewidth=0.8pt]{wl}{zl}
   \ncline[linewidth=0.8pt]{yl}{zl}
   \ncline[linewidth=0.8pt]{yl}{zlb}
    \ncline[linewidth=0.8pt]{vl}{zlb}
    \ncline[linewidth=0.8pt]{vl}{wl}
    \rput(2.05,-1.5){\scriptsize(a)}
 \end{pspicture}
 \begin{pspicture}(-5.6,1)(-4.9,2.2)

   \cnode(0.95,2.05){.09}{xl}\rput(0.55,2.2){\scriptsize$0011$}
  \cnode(0.95,-.15){.09}{vl}\rput(0.55,-.3){\scriptsize$0100$}
  \cnode(3.15,2.05){.09}{ul}\rput(3.55,2.2){\scriptsize$0110$}
  \cnode(3.15,-.15){.09}{yl}\rput(3.55,-0.3){\scriptsize$0001$}
  \cnode(0.5,0.95){.09}{wl} \rput(0.,0.95){\scriptsize$0101$}
  \cnode(3.6,0.95){.09}{zl} \rput(4.1,0.95){\scriptsize$0111$}
  \cnode(2.05,2.5){.09}{wla}\rput(2.05,2.75){\scriptsize$0010$}
  \cnode(2.05,-0.6){.09}{zlb}\rput(2.05,-0.85){\scriptsize$0000$}

 \cnode(5.95,2.05){.09}{rxl}\rput(5.55,2.2){\scriptsize$1110$}
  \cnode(5.95,-.15){.09}{rvl}\rput(5.55,-.3){\scriptsize$1001$}
  \cnode(8.15,2.05){.09}{rul}\rput(8.55,2.2){\scriptsize$1011$}
  \cnode(8.15,-.15){.09}{ryl}\rput(8.55,-0.3){\scriptsize$1100$}
  \cnode(5.5,0.95){.09}{rwl}\rput(5.05,1.){\scriptsize$1111$}
  \cnode(8.6,0.95){.09}{rzl}\rput(9.05,0.95){\scriptsize$1101$}
  \cnode(7.05,2.5){.09}{rwla}\rput(7.05,2.75){\scriptsize$1010$}
  \cnode(7.05,-0.6){.09}{rzlb}\rput(7.05,-0.85){\scriptsize$1000$}

  \ncline[linewidth=0.8pt]{xl}{wla}
  \ncline[linewidth=0.8pt]{wl}{xl}
   \ncline[linewidth=0.8pt]{yl}{xl}
   \ncline[linewidth=0.8pt]{zlb}{wla}
   \ncline[linewidth=0.8pt]{ul}{wla}
   \ncline[linewidth=0.8pt]{ul}{zl}
   \ncline[linewidth=0.8pt]{ul}{vl}
   \ncline[linewidth=0.8pt]{wl}{zl}
   \ncline[linewidth=0.8pt]{yl}{zl}
   \ncline[linewidth=0.8pt]{yl}{zlb}
    \ncline[linewidth=0.8pt]{vl}{zlb}
    \ncline[linewidth=0.8pt]{vl}{wl}

    \ncline[linewidth=0.8pt]{rxl}{rwla}
  \ncline[linewidth=0.8pt]{rwl}{rxl}
   \ncline[linewidth=0.8pt]{ryl}{rxl}
   \ncline[linewidth=0.8pt]{rzlb}{rwla}
   \ncline[linewidth=0.8pt]{rul}{rwla}
   \ncline[linewidth=0.8pt]{rul}{rzl}
   \ncline[linewidth=0.8pt]{rul}{rvl}
   \ncline[linewidth=0.8pt]{rwl}{rzl}
   \ncline[linewidth=0.8pt]{ryl}{rzl}
   \ncline[linewidth=0.8pt]{ryl}{rzlb}
    \ncline[linewidth=0.8pt]{rvl}{rzlb}
    \ncline[linewidth=0.8pt]{rvl}{rwl}

 \ncline[linewidth=.4pt]{wla}{rwla}
  \ncline[linewidth=.4pt]{ul}{rxl}
   \nccurve[angleA=21,angleB=158,linewidth=.4pt]{zl}{rzl}
   \ncline[linewidth=.4pt]{yl}{rul}
   \ncline[linewidth=.4pt]{zlb}{rzlb}
   \nccurve[angleA=21,angleB=158,linewidth=.4pt]{vl}{ryl}
   \nccurve[angleA=21,angleB=158,linewidth=.4pt]{wl}{rwl}
   \ncline[linewidth=.4pt]{xl}{rvl}
 \rput(4.05,-1.5){\scriptsize (b)}
 \end{pspicture}
 \vskip50pt
\caption{
\label{f4}                                  
\footnotesize (a) $CQ_3$;\ \ (b) $CQ_4$}
\end{figure}
\vskip10pt
\begin{figure}[ht]
\psset{unit=0.8}
\begin{pspicture}(-2,1)(0.5,2.2)

  \cnode(0.95,2.05){.09}{xl}\rput(0.6,2.2){\scriptsize$001$}
  \cnode(0.95,-.15){.09}{vl}\rput(0.6,-.3){\scriptsize$110$}
  \cnode(3.15,2.05){.09}{ul}\rput(3.5,2.2){\scriptsize$100$}
  \cnode(3.15,-.15){.09}{yl}\rput(3.5,-0.3){\scriptsize$011$}
  \cnode(0.5,0.95){.09}{wl} \rput(0.15,0.95){\scriptsize$111$}
  \cnode(3.6,0.95){.09}{zl} \rput(3.95,0.95){\scriptsize$101$}
  \cnode(2.05,2.5){.09}{wla}\rput(2.05,2.75){\scriptsize$000$}
  \cnode(2.05,-0.6){.09}{zlb}\rput(2.05,-0.85){\scriptsize$010$}

  \ncline[linewidth=0.8pt]{xl}{wla}
  \ncline[linewidth=0.8pt]{wl}{xl}
   \ncline[linewidth=0.8pt]{yl}{xl}
   \ncline[linewidth=0.8pt]{zlb}{wla}
   \ncline[linewidth=0.8pt]{ul}{wla}
   \ncline[linewidth=0.8pt]{ul}{zl}
   \ncline[linewidth=0.8pt]{ul}{vl}
   \ncline[linewidth=0.8pt]{wl}{zl}
   \ncline[linewidth=0.8pt]{yl}{zl}
   \ncline[linewidth=0.8pt]{yl}{zlb}
    \ncline[linewidth=0.8pt]{vl}{zlb}
    \ncline[linewidth=0.8pt]{vl}{wl}
    \rput(2.05,-1.5){\scriptsize(a)}
 \end{pspicture}
 \begin{pspicture}(-5.6,1)(-4.9,2.2)

   \cnode(0.95,2.05){.09}{xl}\rput(0.55,2.2){\scriptsize$0001$}
  \cnode(0.95,-.15){.09}{vl}\rput(0.55,-.3){\scriptsize$0110$}
  \cnode(3.15,2.05){.09}{ul}\rput(3.55,2.2){\scriptsize$0100$}
  \cnode(3.15,-.15){.09}{yl}\rput(3.55,-0.3){\scriptsize$0011$}
  \cnode(0.5,0.95){.09}{wl} \rput(0.,0.95){\scriptsize$0111$}
  \cnode(3.6,0.95){.09}{zl} \rput(4.1,0.88){\scriptsize$0101$}
  \cnode(2.05,2.5){.09}{wla}\rput(2.05,2.75){\scriptsize$0000$}
  \cnode(2.05,-0.6){.09}{b}\rput(2.05,-0.85){\scriptsize$0010$}

 \cnode(5.95,2.05){.09}{rxl}\rput(5.55,2.2){\scriptsize$1100$}
  \cnode(5.95,-.15){.09}{rvl}\rput(5.45,-.1){\scriptsize$1011$}
  \cnode(8.15,2.05){.09}{rul}\rput(8.55,2.2){\scriptsize$1001$}
  \cnode(8.15,-.15){.09}{ryl}\rput(8.55,-0.3){\scriptsize$1110$}
  \cnode(5.5,0.95){.09}{rwl}\rput(5.05,0.87){\scriptsize$1101$}
  \cnode(8.6,0.95){.09}{rzl}\rput(9.05,0.95){\scriptsize$1111$}
  \cnode(7.05,2.5){.09}{rwla}\rput(7.05,2.75){\scriptsize$1000$}
  \cnode(7.05,-0.6){.09}{rb}\rput(7.05,-0.85){\scriptsize$1010$}

   \ncline[linewidth=0.8pt]{xl}{wla}
  \ncline[linewidth=0.8pt]{wl}{xl}
   \ncline[linewidth=0.8pt]{yl}{xl}
   \ncline[linewidth=0.8pt]{b}{wla}
   \ncline[linewidth=0.8pt]{ul}{wla}
   \ncline[linewidth=0.8pt]{ul}{zl}
   \ncline[linewidth=0.8pt]{ul}{vl}
   \ncline[linewidth=0.8pt]{wl}{zl}
   \ncline[linewidth=0.8pt]{yl}{zl}
   \ncline[linewidth=0.8pt]{yl}{b}
    \ncline[linewidth=0.8pt]{vl}{b}
    \ncline[linewidth=0.8pt]{vl}{wl}

    \ncline[linewidth=0.8pt]{rxl}{rwla}
  \ncline[linewidth=0.8pt]{rwl}{rxl}
   \ncline[linewidth=0.8pt]{ryl}{rxl}
   \ncline[linewidth=0.8pt]{rb}{rwla}
   \ncline[linewidth=0.8pt]{rul}{rwla}
   \ncline[linewidth=0.8pt]{rul}{rzl}
   \ncline[linewidth=0.8pt]{rul}{rvl}
   \ncline[linewidth=0.8pt]{rwl}{rzl}
   \ncline[linewidth=0.8pt]{ryl}{rzl}
   \ncline[linewidth=0.8pt]{ryl}{rb}
    \ncline[linewidth=0.8pt]{rvl}{rb}
    \ncline[linewidth=0.8pt]{rvl}{rwl}

 \ncline[linewidth=.4pt]{wla}{rwla}
  \ncline[linewidth=.4pt]{ul}{rxl}
   \ncline[linewidth=.4pt]{zl}{rul}
   \ncline[linewidth=.4pt]{yl}{rzl}
   \ncline[linewidth=.4pt]{b}{rb}
   \nccurve[angleA=21,angleB=158,linewidth=.4pt]{vl}{ryl}
  \ncline[linewidth=.4pt]{wl}{rvl}
  \ncline[linewidth=.4pt]{xl}{rwl}
 \rput(4.05,-1.5){\scriptsize (b)}
 \end{pspicture}
 \vskip50pt
\caption{
\label{f4}                                    
\footnotesize (a) $LTQ_3$;\ \ (b) $LTQ_4$}
\end{figure}
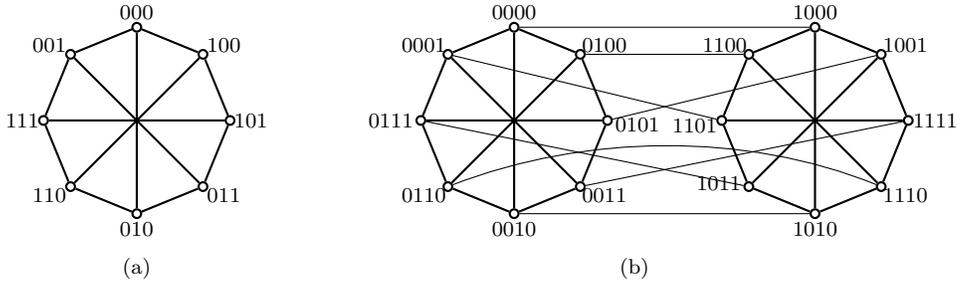

\begin{cor}\textnormal{If $F\subset V(MQ_n)\cup E(MQ_n)$ with $|F|\leq n-2$, then for any two different correct vertices $u$ and $v$, there exists a fault-free path $P_{uv}$ of every length $l$ with $2^{n-1}-1\leq l\leq 2^n-f_v-1-\alpha$, where $\alpha=0$ if vertices $u$ and $v$ form a normal vertex-pair and $\alpha=1$ if vertices $u$ and $v$ form a weak vertex-pair in $MQ_n-F$($n\geq5$). }
\end{cor}
\begin{figure}[ht]
\psset{unit=0.8}
\begin{pspicture}(-2,1)(0.5,2.2)

  \cnode(0.95,2.05){.09}{xl}\rput(0.6,2.2){\scriptsize$000$}
  \cnode(0.95,-.15){.09}{vl}\rput(0.6,-.3){\scriptsize$111$}
  \cnode(3.15,2.05){.09}{ul}\rput(3.5,2.2){\scriptsize$110$}
  \cnode(3.15,-.15){.09}{yl}\rput(3.5,-0.3){\scriptsize$001$}
  \cnode(0.5,0.95){.09}{wl} \rput(0.15,0.95){\scriptsize$100$}
  \cnode(3.6,0.95){.09}{zl} \rput(3.95,0.95){\scriptsize$101$}
  \cnode(2.05,2.5){.09}{wla}\rput(2.05,2.75){\scriptsize$010$}
  \cnode(2.05,-0.6){.09}{zlb}\rput(2.05,-0.85){\scriptsize$011$}

  \ncline[linewidth=0.8pt]{xl}{wla}
  \ncline[linewidth=0.8pt]{wl}{xl}
   \ncline[linewidth=0.8pt]{yl}{xl}
   \ncline[linewidth=0.8pt]{zlb}{wla}
   \ncline[linewidth=0.8pt]{ul}{wla}
   \ncline[linewidth=0.8pt]{ul}{zl}
   \ncline[linewidth=0.8pt]{ul}{vl}
   \ncline[linewidth=0.8pt]{wl}{zl}
   \ncline[linewidth=0.8pt]{yl}{zl}
   \ncline[linewidth=0.8pt]{yl}{zlb}
    \ncline[linewidth=0.8pt]{vl}{zlb}
    \ncline[linewidth=0.8pt]{vl}{wl}
    \rput(2.05,-1.5){\scriptsize(a)}
 \end{pspicture}
 \begin{pspicture}(-5.6,1)(-4.9,2.2)

   \cnode(0.95,2.05){.09}{xl}\rput(0.55,2.2){\scriptsize$0000$}
  \cnode(0.95,-.15){.09}{vl}\rput(0.55,-.3){\scriptsize$0111$}
  \cnode(3.15,2.05){.09}{ul}\rput(3.55,2.2){\scriptsize$0110$}
  \cnode(3.15,-.15){.09}{yl}\rput(3.55,-0.3){\scriptsize$0001$}
  \cnode(0.5,0.95){.09}{wl} \rput(0.,0.95){\scriptsize$0100$}
  \cnode(3.6,0.95){.09}{zl} \rput(4.1,0.95){\scriptsize$0101$}
  \cnode(2.05,2.5){.09}{wla}\rput(2.05,2.75){\scriptsize$0010$}
  \cnode(2.05,-0.6){.09}{b}\rput(2.05,-0.85){\scriptsize$0011$}

 \cnode(5.95,2.05){.09}{rxl}\rput(5.55,2.2){\scriptsize$1101$}
  \cnode(5.95,-.15){.09}{rvl}\rput(5.55,-.3){\scriptsize$1001$}
  \cnode(8.15,2.05){.09}{rul}\rput(8.55,2.2){\scriptsize$1000$}
  \cnode(8.15,-.15){.09}{ryl}\rput(8.55,-0.3){\scriptsize$1100$}
  \cnode(5.5,0.95){.09}{rwl}\rput(5.05,1.){\scriptsize$1110$}
  \cnode(8.6,0.95){.09}{rzl}\rput(9.05,0.95){\scriptsize$1111$}
  \cnode(7.05,2.5){.09}{rwla}\rput(7.05,2.75){\scriptsize$1010$}
  \cnode(7.05,-0.6){.09}{rb}\rput(7.05,-0.85){\scriptsize$1011$}

  \ncline[linewidth=0.8pt]{xl}{wla}
  \ncline[linewidth=0.8pt]{wl}{xl}
   \ncline[linewidth=0.8pt]{yl}{xl}
   \ncline[linewidth=0.8pt]{b}{wla}
   \ncline[linewidth=0.8pt]{ul}{wla}
   \ncline[linewidth=0.8pt]{ul}{zl}
   \ncline[linewidth=0.8pt]{ul}{vl}
   \ncline[linewidth=0.8pt]{wl}{zl}
   \ncline[linewidth=0.8pt]{yl}{zl}
   \ncline[linewidth=0.8pt]{yl}{b}
    \ncline[linewidth=0.8pt]{vl}{b}
    \ncline[linewidth=0.8pt]{vl}{wl}

    \ncline[linewidth=0.8pt]{rxl}{rwla}
  \ncline[linewidth=0.8pt]{rwl}{rxl}
   \ncline[linewidth=0.8pt]{ryl}{rxl}
   \ncline[linewidth=0.8pt]{rb}{rwla}
   \ncline[linewidth=0.8pt]{rul}{rwla}
   \ncline[linewidth=0.8pt]{rul}{rzl}
   \ncline[linewidth=0.8pt]{rul}{rvl}
   \ncline[linewidth=0.8pt]{rwl}{rzl}
   \ncline[linewidth=0.8pt]{ryl}{rzl}
   \ncline[linewidth=0.8pt]{ryl}{rb}
    \ncline[linewidth=0.8pt]{rvl}{rb}
    \ncline[linewidth=0.8pt]{rvl}{rwl}

 \ncline[linewidth=.4pt]{wla}{rwla}
  \ncline[linewidth=.4pt]{ul}{rwl}
   \ncline[linewidth=.4pt]{zl}{rxl}
   \ncline[linewidth=.4pt]{yl}{rvl}
   \ncline[linewidth=.4pt]{b}{rb}
   \ncline[linewidth=.4pt]{vl}{rzl}
  \ncline[linewidth=.4pt]{wl}{ryl}
  \nccurve[angleA=-21,angleB=-158,linewidth=.4pt]{xl}{rul}
 \rput(4.05,-1.5){\scriptsize (b)}
 \end{pspicture}
 \vskip80pt

\begin{pspicture}(-2,1)(0.5,2.2)

  \cnode(0.95,2.05){.09}{xl}\rput(0.6,2.2){\scriptsize$000$}
  \cnode(0.95,-.15){.09}{vl}\rput(0.6,-.3){\scriptsize$111$}
  \cnode(3.15,2.05){.09}{ul}\rput(3.5,2.2){\scriptsize$110$}
  \cnode(3.15,-.15){.09}{yl}\rput(3.5,-0.3){\scriptsize$001$}
  \cnode(0.5,0.95){.09}{wl} \rput(0.15,0.95){\scriptsize$100$}
  \cnode(3.6,0.95){.09}{zl} \rput(3.95,0.95){\scriptsize$101$}
  \cnode(2.05,2.5){.09}{wla}\rput(2.05,2.75){\scriptsize$010$}
  \cnode(2.05,-0.6){.09}{zlb}\rput(2.05,-0.85){\scriptsize$011$}

  \ncline[linewidth=0.8pt]{xl}{wla}
  \ncline[linewidth=0.8pt]{wl}{xl}
   \ncline[linewidth=0.8pt]{yl}{xl}
   \ncline[linewidth=0.8pt]{zlb}{wla}
   \ncline[linewidth=0.8pt]{ul}{wla}
   \ncline[linewidth=0.8pt]{ul}{zl}
   \ncline[linewidth=0.8pt]{ul}{vl}
   \ncline[linewidth=0.8pt]{wl}{zl}
   \ncline[linewidth=0.8pt]{yl}{zl}
   \ncline[linewidth=0.8pt]{yl}{zlb}
    \ncline[linewidth=0.8pt]{vl}{zlb}
    \ncline[linewidth=0.8pt]{vl}{wl}
    \rput(2.05,-1.5){\scriptsize(c)}
 \end{pspicture}
 \begin{pspicture}(-5.6,1)(-4.9,2.2)

   \cnode(0.95,2.05){.09}{xl}\rput(0.55,2.2){\scriptsize$0000$}
  \cnode(0.95,-.15){.09}{vl}\rput(0.55,-.3){\scriptsize$0111$}
  \cnode(3.15,2.05){.09}{ul}\rput(3.55,2.2){\scriptsize$0110$}
  \cnode(3.15,-.15){.09}{yl}\rput(3.55,-0.3){\scriptsize$0001$}
  \cnode(0.5,0.95){.09}{wl} \rput(0.,0.95){\scriptsize$0100$}
  \cnode(3.6,0.95){.09}{zl} \rput(4.1,0.95){\scriptsize$0101$}
  \cnode(2.05,2.5){.09}{wla}\rput(2.05,2.75){\scriptsize$0010$}
  \cnode(2.05,-0.6){.09}{b}\rput(2.05,-0.85){\scriptsize$0011$}

 \cnode(5.95,2.15){.09}{rxl}\rput(5.55,2.3){\scriptsize$1101$}
  \cnode(5.95,-.15){.09}{rvl}\rput(5.45,-.1){\scriptsize$1001$}
  \cnode(8.15,2.05){.09}{rul}\rput(8.55,2.2){\scriptsize$1000$}
  \cnode(8.15,-.15){.09}{ryl}\rput(8.55,-0.3){\scriptsize$1100$}
  \cnode(5.5,0.95){.09}{rwl}\rput(5.05,1.){\scriptsize$1110$}
  \cnode(8.6,0.95){.09}{rzl}\rput(9.05,0.95){\scriptsize$1111$}
  \cnode(7.05,2.5){.09}{rwla}\rput(7.05,2.75){\scriptsize$1010$}
  \cnode(7.05,-0.6){.09}{rb}\rput(7.05,-0.85){\scriptsize$1011$}

   \ncline[linewidth=0.8pt]{xl}{wla}
  \ncline[linewidth=0.8pt]{wl}{xl}
   \ncline[linewidth=0.8pt]{yl}{xl}
   \ncline[linewidth=0.8pt]{b}{wla}
   \ncline[linewidth=0.8pt]{ul}{wla}
   \ncline[linewidth=0.8pt]{ul}{zl}
   \ncline[linewidth=0.8pt]{ul}{vl}
   \ncline[linewidth=0.8pt]{wl}{zl}
   \ncline[linewidth=0.8pt]{yl}{zl}
   \ncline[linewidth=0.8pt]{yl}{b}
    \ncline[linewidth=0.8pt]{vl}{b}
    \ncline[linewidth=0.8pt]{vl}{wl}

    \ncline[linewidth=0.8pt]{rxl}{rwla}
  \ncline[linewidth=0.8pt]{rwl}{rxl}
   \ncline[linewidth=0.8pt]{ryl}{rxl}
   \ncline[linewidth=0.8pt]{rb}{rwla}
   \ncline[linewidth=0.8pt]{rul}{rwla}
   \ncline[linewidth=0.8pt]{rul}{rzl}
   \ncline[linewidth=0.8pt]{rul}{rvl}
   \ncline[linewidth=0.8pt]{rwl}{rzl}
   \ncline[linewidth=0.8pt]{ryl}{rzl}
   \ncline[linewidth=0.8pt]{ryl}{rb}
    \ncline[linewidth=0.8pt]{rvl}{rb}
    \ncline[linewidth=0.8pt]{rvl}{rwl}

 \ncline[linewidth=.4pt]{wla}{rxl}
 \ncline[linewidth=.4pt]{ul}{rvl}
  \ncline[linewidth=.4pt]{zl}{rwla}
   \ncline[linewidth=.4pt]{yl}{rwl}
   \ncline[linewidth=.4pt]{b}{ryl}
   \ncline[linewidth=.4pt]{vl}{rul}
   \ncline[linewidth=.4pt]{wl}{rb}
  \ncline[linewidth=.4pt]{xl}{rzl}
 \rput(4.05,-1.5){\scriptsize (d)}
 \end{pspicture}
 \vskip50pt
\caption{
\label{f4}                                
\footnotesize (a) $MQ^0_3$;\ \ (b) $MQ^0_4$;\ \ (c) $MQ^1_3$;\ \ (d) $MQ^1_4$}
\end{figure}
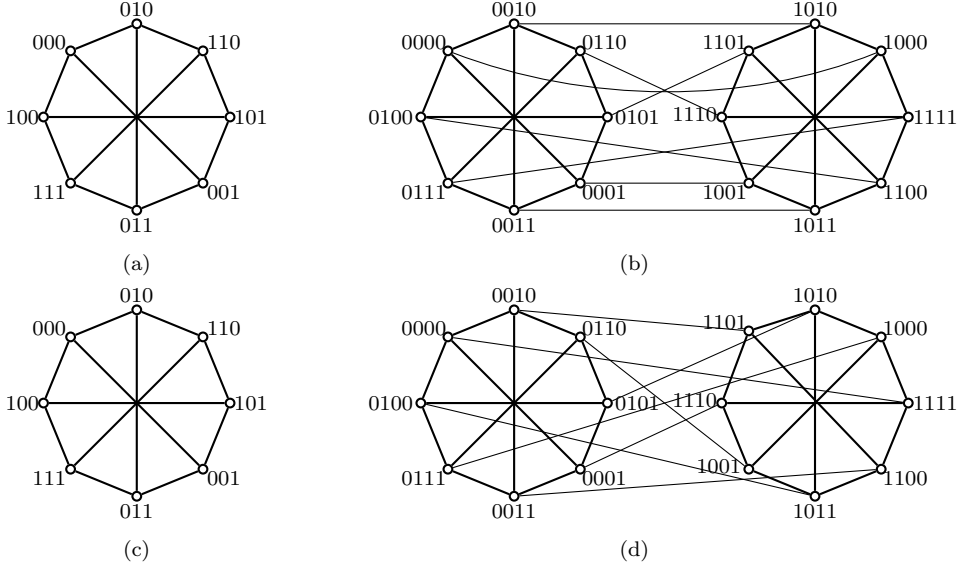
\begin{cor}\textnormal{If $F\subset V(TQ_n)\cup E(TQ_n)$ with $|F|\leq n-2$, then for any two different correct vertices $u$ and $v$, there exists a fault-free path $P_{uv}$ of every length $l$ with $2^{n-1}-1\leq l\leq 2^n-f_v-1-\alpha$, where $\alpha=0$ if vertices $u$ and $v$ form a normal vertex-pair and $\alpha=1$ if vertices $u$ and $v$ form a weak vertex-pair in $TQ_n-F$($n\geq5$) for any odd $n$. }
\end{cor}
In this paper, we apply our strategy to these four network topologies($MQ_n,LTQ_n,\\TQ_n,CQ_n$). In the future work, we will extend our strategy to other
graphs of Hypercube-Like Networks.
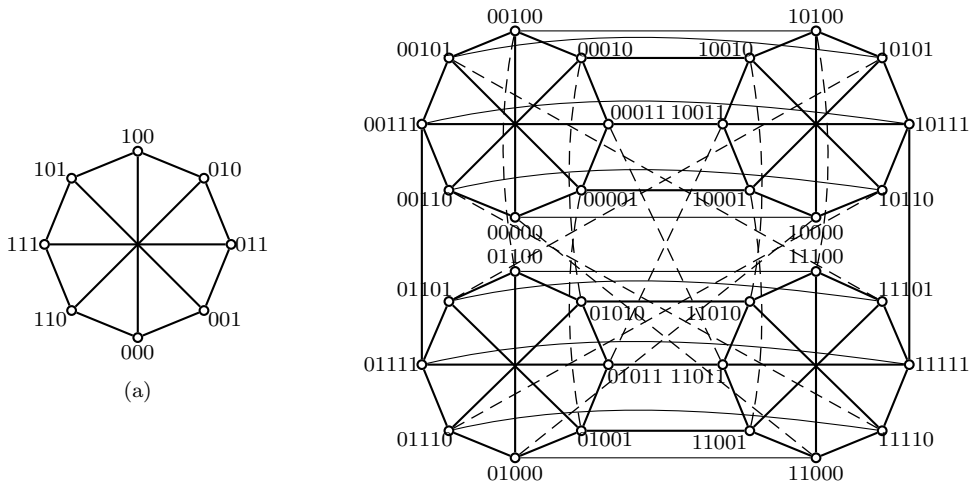
\begin{figure}[ht]
\psset{unit=0.8}
\begin{pspicture}(-2,3)(0.5,4.2)

  \cnode(0.95,2.05){.09}{xl}\rput(0.6,2.2){\scriptsize$101$}
  \cnode(0.95,-.15){.09}{vl}\rput(0.6,-.3){\scriptsize$110$}
  \cnode(3.15,2.05){.09}{ul}\rput(3.5,2.2){\scriptsize$010$}
  \cnode(3.15,-.15){.09}{yl}\rput(3.5,-0.3){\scriptsize$001$}
  \cnode(0.5,0.95){.09}{wl} \rput(0.15,0.95){\scriptsize$111$}
  \cnode(3.6,0.95){.09}{zl} \rput(3.95,0.95){\scriptsize$011$}
  \cnode(2.05,2.5){.09}{wla}\rput(2.05,2.75){\scriptsize$100$}
  \cnode(2.05,-0.6){.09}{zlb}\rput(2.05,-0.85){\scriptsize$000$}

  \ncline[linewidth=0.8pt]{xl}{wla}
  \ncline[linewidth=0.8pt]{wl}{xl}
   \ncline[linewidth=0.8pt]{yl}{xl}
   \ncline[linewidth=0.8pt]{zlb}{wla}
   \ncline[linewidth=0.8pt]{ul}{wla}
   \ncline[linewidth=0.8pt]{ul}{zl}
   \ncline[linewidth=0.8pt]{ul}{vl}
   \ncline[linewidth=0.8pt]{wl}{zl}
   \ncline[linewidth=0.8pt]{yl}{zl}
   \ncline[linewidth=0.8pt]{yl}{zlb}
    \ncline[linewidth=0.8pt]{vl}{zlb}
    \ncline[linewidth=0.8pt]{vl}{wl}
    \rput(2.05,-1.5){\scriptsize(a)}
 \end{pspicture}
 \begin{pspicture}(-5.6,1)(-4.9,2.2)

   \cnode(0.95,2.05){.09}{xl}\rput(0.55,2.2){\scriptsize$00101$}
  \cnode(0.95,-.15){.09}{vl}\rput(0.55,-.3){\scriptsize$00110$}
  \cnode(3.15,2.05){.09}{ul}\rput(3.55,2.2){\scriptsize$00010$}
  \cnode(3.15,-.15){.09}{yl}\rput(3.65,-0.3){\scriptsize$00001$}
  \cnode(0.5,0.95){.09}{wl} \rput(0.,0.95){\scriptsize$00111$}
  \cnode(3.6,0.95){.09}{zl} \rput(4.1,1.15){\scriptsize$00011$}
  \cnode(2.05,2.5){.09}{wla}\rput(2.05,2.75){\scriptsize$00100$}
  \cnode(2.05,-0.6){.09}{b}\rput(2.05,-0.85){\scriptsize$00000$}

 \cnode(5.95,2.05){.09}{ rul}\rput(5.55,2.2){\scriptsize$10010$}
  \cnode(5.95,-.15){.09}{ryl}\rput(5.45,-.3){\scriptsize$10001$}
  \cnode(8.15,2.05){.09}{rxl}\rput(8.55,2.2){\scriptsize$10101$}
  \cnode(8.15,-.15){.09}{rvl}\rput(8.55,-0.3){\scriptsize$10110$}
  \cnode(5.5,0.95){.09}{rzl}\rput(5.1,1.15){\scriptsize$10011$}
  \cnode(8.6,0.95){.09}{rwl}\rput(9.15,0.95){\scriptsize$10111$}
  \cnode(7.05,2.5){.09}{rwla}\rput(7.05,2.75){\scriptsize$10100$}
  \cnode(7.05,-0.6){.09}{rb}\rput(7.05,-0.85){\scriptsize$10000$}

   \ncline[linewidth=0.8pt]{xl}{wla}
  \ncline[linewidth=0.8pt]{wl}{xl}
   \ncline[linewidth=0.8pt]{yl}{xl}
   \ncline[linewidth=0.8pt]{b}{wla}
   \ncline[linewidth=0.8pt]{ul}{wla}
   \ncline[linewidth=0.8pt]{ul}{zl}
   \ncline[linewidth=0.8pt]{ul}{vl}
   \ncline[linewidth=0.8pt]{wl}{zl}
   \ncline[linewidth=0.8pt]{yl}{zl}
   \ncline[linewidth=0.8pt]{yl}{b}
    \ncline[linewidth=0.8pt]{vl}{b}
    \ncline[linewidth=0.8pt]{vl}{wl}

    \ncline[linewidth=0.8pt]{rxl}{rwla}
  \ncline[linewidth=0.8pt]{rwl}{rxl}
   \ncline[linewidth=0.8pt]{ryl}{rxl}
   \ncline[linewidth=0.8pt]{rb}{rwla}
   \ncline[linewidth=0.8pt]{rul}{rwla}
   \ncline[linewidth=0.8pt]{rul}{rzl}
   \ncline[linewidth=0.8pt]{rul}{rvl}
   \ncline[linewidth=0.8pt]{rwl}{rzl}
   \ncline[linewidth=0.8pt]{ryl}{rzl}
   \ncline[linewidth=0.8pt]{ryl}{rb}
    \ncline[linewidth=0.8pt]{rvl}{rb}
    \ncline[linewidth=0.8pt]{rvl}{rwl}

  \ncline[linewidth=.4pt]{wla}{rwla}
  \ncline[linewidth=.4pt]{b}{rb}
  \nccurve[angleA=13,angleB=172,linewidth=.4pt]{xl}{rxl}
  \nccurve[angleA=13,angleB=172,linewidth=.4pt]{vl}{rvl}
   \ncline[linewidth=0.8pt]{ul}{rul}
    \ncline[linewidth=0.8pt]{yl}{ryl}
   \nccurve[angleA=13,angleB=172,linewidth=.4pt]{wl}{rwl}
   \ncline[linewidth=0.8pt]{zl}{rzl}

  \cnode(0.95,-2){.09}{bxl}\rput(0.55,-1.8){\scriptsize$01101$}
  \cnode(0.95,-4.15){.09}{bvl}\rput(0.55,-4.3){\scriptsize$01110$}
  \cnode(3.15,-2.){.09}{bul}\rput(3.75,-2.2){\scriptsize$01010$}
  \cnode(3.15,-4.15){.09}{byl}\rput(3.55,-4.3){\scriptsize$01001$}
  \cnode(0.5,-3.05){.09}{bwl} \rput(0.,-3.05){\scriptsize$01111$}
  \cnode(3.6,-3.05){.09}{bzl} \rput(4.05,-3.25){\scriptsize$01011$}
  \cnode(2.05,-1.5){.09}{bwla}\rput(2.05,-1.25){\scriptsize$01100$}
  \cnode(2.05,-4.6){.09}{bb}\rput(2.05,-4.85){\scriptsize$01000$}

  \cnode(5.95,-2){.09}{brul}\rput(5.35,-2.2){\scriptsize$11010$}
  \cnode(5.95,-4.15){.09}{bryl}\rput(5.45,-4.35){\scriptsize$11001$}
  \cnode(8.15,-2){.09}{brxl}\rput(8.55,-1.8){\scriptsize$11101$}
  \cnode(8.15,-4.15){.09}{brvl}\rput(8.55,-4.3){\scriptsize$11110$}
  \cnode(5.5,-3.05){.09}{brzl}\rput(5.1,-3.25){\scriptsize$11011$}
  \cnode(8.6,-3.05){.09}{brwl}\rput(9.15,-3.05){\scriptsize$11111$}
  \cnode(7.05,-1.5){.09}{brwla}\rput(7.05,-1.25){\scriptsize$11100$}
  \cnode(7.05,-4.6){.09}{brb}\rput(7.05,-4.85){\scriptsize$11000$}

   \ncline[linewidth=0.8pt]{bxl}{bwla}
  \ncline[linewidth=0.8pt]{bwl}{bxl}
   \ncline[linewidth=0.8pt]{byl}{bxl}
   \ncline[linewidth=0.8pt]{bb}{bwla}
   \ncline[linewidth=0.8pt]{bul}{bwla}
   \ncline[linewidth=0.8pt]{bul}{bzl}
   \ncline[linewidth=0.8pt]{bul}{bvl}
   \ncline[linewidth=0.8pt]{bwl}{bzl}
   \ncline[linewidth=0.8pt]{byl}{bzl}
   \ncline[linewidth=0.8pt]{byl}{bb}
    \ncline[linewidth=0.8pt]{bvl}{bb}
    \ncline[linewidth=0.8pt]{bvl}{bwl}

    \ncline[linewidth=0.8pt]{brxl}{brwla}
  \ncline[linewidth=0.8pt]{brwl}{brxl}
   \ncline[linewidth=0.8pt]{bryl}{brxl}
   \ncline[linewidth=0.8pt]{brb}{brwla}
   \ncline[linewidth=0.8pt]{brul}{brwla}
   \ncline[linewidth=0.8pt]{brul}{brzl}
   \ncline[linewidth=0.8pt]{brul}{brvl}
   \ncline[linewidth=0.8pt]{brwl}{brzl}
   \ncline[linewidth=0.8pt]{bryl}{brzl}
   \ncline[linewidth=0.8pt]{bryl}{brb}
    \ncline[linewidth=0.8pt]{brvl}{brb}
    \ncline[linewidth=0.8pt]{brvl}{brwl}

  \ncline[linewidth=.4pt]{bwla}{brwla}
  \ncline[linewidth=.4pt]{bb}{brb}
  \nccurve[angleA=13,angleB=172,linewidth=.4pt]{bxl}{brxl}
  \nccurve[angleA=13,angleB=172,linewidth=.4pt]{bvl}{brvl}
    \ncline[linewidth=0.8pt]{bul}{brul}
   \ncline[linewidth=0.8pt]{byl}{bryl}
   \nccurve[angleA=13,angleB=172,linewidth=.4pt]{bwl}{brwl}
    \ncline[linewidth=0.8pt]{bzl}{brzl}

   \ncline[linewidth=0.8pt]{wl}{bwl}
   \ncline[linewidth=0.8pt]{rwl}{brwl}
   \ncline[linewidth=.5pt, linestyle=dashed]{xl}{brxl}
   \ncline[linewidth=.5pt, linestyle=dashed]{zl}{brzl}
    \ncline[linewidth=.5pt, linestyle=dashed]{b}{brb}
    \ncline[linewidth=.5pt, linestyle=dashed]{vl}{brvl}
    \ncline[linewidth=.5pt, linestyle=dashed]{rxl}{bxl}
     \ncline[linewidth=.5pt, linestyle=dashed]{rvl}{bvl}
      \ncline[linewidth=.5pt, linestyle=dashed]{rzl}{bzl}
      \ncline[linewidth=.5pt, linestyle=dashed]{rb}{bb}
       \nccurve[angleA=-102,angleB=100,linewidth=.5pt, linestyle=dashed]{wla}{bwla}
        \nccurve[angleA=-102,angleB=100,linewidth=.5pt, linestyle=dashed]{yl}{byl}
        \nccurve[angleA=-102,angleB=100,linewidth=.5pt, linestyle=dashed]{ul}{bul}

         \nccurve[angleA=-80,angleB=78,linewidth=.5pt, linestyle=dashed]{rwla}{brwla}
        \nccurve[angleA=-80,angleB=78,linewidth=.5pt, linestyle=dashed]{ryl}{bryl}
        \nccurve[angleA=-80,angleB=78,linewidth=.5pt, linestyle=dashed]{rul}{brul}
 \end{pspicture}
 \vskip130pt
\caption{
\label{f4}                                    
\footnotesize (a) $TQ_3$;\ \ (b) $TQ_5$}
\end{figure}

\end{document}